\documentclass[a4paper]{amsart}
\usepackage{amsmath,amssymb,amsfonts}
\usepackage{colordvi}
\usepackage[all]{xy}

\def\mpoint{\;.}
\def\mvirg{\;,}

\def\mpn{\medskip\par\noindent}

\def\mmpn{\vskip 1em minus 1em\par\noindent}

\def\sp{\bigskip\par}

\def\smp{\smallskip\par}

\def\CB{{\mathcal B}}
\def\CC{{\mathcal C}}
\def\CD{{\mathcal D}}
\def\CE{{\mathcal E}}
\def\CF{{\mathcal F}}

\def\CI{{\mathcal I}}
\def\CL{{\mathcal L}}

\def\CP{{\mathcal P}}

\def\CX{{\mathcal X}}

\def\Ker{\operatorname{Ker}\nolimits}

\def\Im{\operatorname{Im}\nolimits}
\def\Id{\operatorname{id}\nolimits}

\def\Mod{\operatorname{Mod}\nolimits}
\def\Hom{\operatorname{Hom}\nolimits}
\def\Inj{\operatorname{Inj}\nolimits}
\def\Sur{\operatorname{Sur}\nolimits}
\def\End{\operatorname{End}\nolimits}
\def\Aut{\operatorname{Aut}\nolimits}

\def\Irr{\operatorname{Irr}\nolimits}

\def\Pol{\operatorname{Pol}\nolimits}

\def\op{^{op}}

\def\ls#1#2{{\,^{#1}\!#2}}
\def\Iup{I^{\uparrow}}
\def\Idown{I_{\downarrow}}
\def\meet{\wedge}

\def\Z{\mathbb{Z}}

\def\F{\mathbb{F}}

\def\S{\mathbb{S}}

\newcommand{\edge}[2]{\xymatrix{#1\ar@{->-}[r]&#2}}
\def\marc[#1]{\ar@{-}[#1]|(.4){\object@{<}}}
\def\mard[#1]{\ar@{-}[#1]|(.5){\object@{>}}}
\def\marb[#1]{\ar@{-}[#1]|{\object+{  }}}
\newcommand{\fleche}[2]{\xymatrix@C=4ex{*!U(0.2){#1\;}&*!U(0.5){\;#2}\marc[l]}}
\newcommand{\flecheb}[2]{\xymatrix@C=4ex{*!U(0.2){#1\;}&*!U(0.1){\;#2}\marc[l]}}

\def\pf{\par\bigskip\noindent{\bf Proof~: }}
\def\endpf{~\hfill\rlap{\hspace{-1ex}\raisebox{.5ex}{\framebox[1ex]{}}\sp}\bigskip\pagebreak[3]}

\makeatletter
\renewenvironment{enumerate}{\ifnum \@enumdepth >3 \@toodeep\else
       \advance\@enumdepth \@ne
       \edef\@enumctr{enum\romannumeral\the\@enumdepth}\list
       {\csname  label\@enumctr\endcsname}{\setlength{\topsep}{1ex}
 \setlength{\itemsep}{0 pt}\usecounter
         {\@enumctr}\def\makelabel##1{\hss\llap{##1}}}\fi}{\endlist}

\def\@seccntformat#1{\csname the#1\endcsname.\quad}

\def\section{\pagebreak[3]\setcounter{prop}{0}\setcounter{equation}{0}\@startsection{section}{1}{\z@}{4ex plus  6ex}{2ex}{\center\reset@font \large\bf}}

\makeatother

\def\theprop{\thesection.\arabic{prop}}

\renewenvironment{equation}{\refstepcounter{subsection}\refstepcounter
{prop}$$}{\leqno{\bf (\theprop)}$$}

\newenvironment{rem}[1]{\refstepcounter{subsection}\refstepcounter
{prop} \mpn{{\bf \thesection.\arabic{prop}.}\ \ \bf#1.}}{\smp}
\newenvironment{enonce}[1]{\pagebreak[3]\refstepcounter{prop}\mmpn
{{\bf  \thesection.\arabic{prop}.\ #1.}}\begin{it} }{\end{it}\smp}
\def\thesection{\arabic{section}}

\newcommand{\result}[1]{\begin{enonce}{#1}}
\newcommand{\fresult}{\end{enonce}}

\begin{document}

\title[Simple and projective correspondence functors]
{Simple and projective correspondence functors}

\author{Serge Bouc}
\author{Jacques Th\'evenaz}
\date\today



\begin{abstract}
A correspondence functor is a functor from the category of finite sets and correspondences
to the category of $k$-modules, where $k$ is a commutative ring.
We determine exactly which simple correspondence functors are projective.
Moreover, we analyze the occurrence of such simple projective functors inside the correspondence functor~$F$ associated with a finite lattice
and we deduce a direct sum decomposition of~$F$.
\end{abstract}

\maketitle

%%%%%%%% Section

\section{Introduction}

\bigskip
\noindent
In the present paper, we continue to develop the theory of correspondence functors,
namely functors from the category of finite sets and correspondences to the category of $k$-modules, where $k$ is a commutative ring.
Assuming that $k$ is a field, we showed in~\cite{BT2} how to parametrize the simple correspondence functors $S_{E,R,V}$
by means of a finite set $E$, an order relation $R$ on~$E$, and a simple $k\Aut(E,R)$-module~$V$ (up to isomorphism).
Here, we determine which of them are projective (or equivalently injective).\par

We say that a poset $(E,R)$ is a {\em pole poset} if it is obtained by stacking posets having either cardinality one or cardinality two with two incomparable elements (see Section~\ref{Section-pole} for details).

\result{Theorem} 
Let $k$ be a field and let $S_{E,R,V}$ be the simple correspondence functor
parametrized by a finite set $E$, an order relation $R$ on~$E$, and a simple $k\Aut(E,R)$-module~$V$.
The following conditions are equivalent~:
\begin{enumerate}
\item $S_{E,R,V}$ is projective.
\item The poset $(E,R)$ is a pole poset and $V$ is a projective $k\Aut(E,R)$-module.
\item Either $(E,R)$ is a totally ordered poset or $(E,R)$ is a pole poset and the characteristic of~$k$ is different from~2.
\end{enumerate}
\fresult

Since the group $\Aut(E,R)$ of automorphisms of a pole poset is a 2-group, (b) and (c) are easily seen to be equivalent.
However, it requires much more work to prove that (a) implies (b), and also that (b) implies (a) (see Section~\ref{Section-characterization}).
In the case when $(E,R)$ is totally ordered, the projectivity of~$S_{E,R,V}$ was already proved in Corollary 11.11 of~\cite{BT3}.\par

Every simple functor $S_{E,R,V}$ has a precursor $\S_{E,R}$,
called the fundamental functor associated with the poset $(E,R)$ (see Proposition~\ref{precursor}).
This functor $\S_{E,R}$ has the advantage of being defined over any commutative base ring~$k$.
In analogy with the theorem above, we prove in Section~\ref{Section-fundamental} that $\S_{E,R}$ is projective if and only if $(E,R)$ is a pole poset.\par

Associated with a finite lattice~$T$, there is a correspondence functor $F_T$ which is defined over an arbitrary commutative ring~$k$
and which plays a crucial role in the theory, see \cite{BT3, BT4}.
We know in particular that $F_T$ is projective if and only if the lattice $T$ is distributive, for instance if $T$ is a pole lattice.
Also, the assignment $T\mapsto F_T$ is known to be a fully faithful functor by~\cite{BT3}.\par

If $T$ is arbitrary, we show that $F_T$ has direct summands corresponding to pole lattices appearing inside~$T$,
by means of suitably constructed idempotents in~$\End(F_T)$.
Actually, most of the work is done in~$\End(T)$ (where morphisms between lattices are defined to be $k$-linear combinations of join-morphisms), and then corresponding results for~$\End(F_T)$ are obtained using the fully faithful functor $T\mapsto F_T$.
The construction of idempotents in~$\End(T)$ is quite technical (see Section~\ref{Section-idempotents})
but it provides an explicit description of the part of~$\End(T)$ which corresponds to pole lattices (see Section~\ref{Section-subalgebras}).\par

In Section~\ref{Section-functors-pole}, we analyze the special case when $Q$ is a pole lattice (see Theorem~\ref{endom-pole} for details).

\result{Theorem} Let $Q$ be a pole lattice. Then $\End(F_Q)$ is isomorphic to a direct sum of matrix algebras
$$\End(F_Q) \cong\End(Q) \cong \bigoplus_{P} M_{n(Q,P)}(k\Aut(P))$$
where $P$ varies among pole lattices inside~$Q$ and $n(Q,P)$ is some explicit integer.
\fresult

From this, we obtain a decomposition of~$F_Q$ as a direct sum of projective functors (Theorem~\ref{functor-SQ}) and each summand is also simple when $k$ is a field
(Corollary~\ref{simple-SQ}).
Finally, if $T$ is an arbitrary finite lattice, we describe a projective direct summand of~$F_T$ corresponding to all pole lattices which appear inside~$T$
(Theorem~\ref{pole-direct-factors}).

%%%%%%%% Section

\section{Pole posets, pole lattices, and opposite morphisms} \label{Section-pole}
\bigskip
\noindent
We first recall some standard facts about lattices and fix the terminology and the notation.
If $T$ is a finite lattice, we denote by $\vee$ its join, $\meet$ its meet, and $\leq_T$ its order relation.
When the context is clear, we simply write $\leq$ instead of~$\leq_T$.
The unique minimal element is written~$\widehat0$ and the unique maximal element~$\widehat1$.
We let $T\op$ denote the opposite lattice, such that
$$x\leq_T y \iff y\leq_{T\op} x \mpoint$$
A join-irreducible element in~$T$ is simply called {\em irreducible}.
We write $\Irr(T)$ for the full subposet of irreducible elements of~$T$.
Recall that $\widehat0$ is an empty join, hence is not irreducible.
Similarly $\widehat1$ is an empty meet.
If $e\in T$ is irreducible, then the half-open interval $[\,\widehat0,e\,[_T$ has a unique maximal element, written~$r(e)$.
In other words, $r(e)=\sup\{x\in T\mid x<e\}$.
Similarly, if $a$ is meet-irreducible (i.e. irreducible in the opposite lattice~$T\op$), then we define $s(a)=\inf\{x\in T\mid a<x\}$.
Any finite poset~$A$ is isomorphic to the full subposet of irreducible elements of a lattice,
e.g. the lattice~$\Idown(A)$ of all subsets of~$A$ closed under taking smaller elements.\par

Now we want to introduce one of the main concepts for the present a paper.
Let $A$ and $B$ be two finite posets. Define $A*B$ to be the poset whose underlying set is the disjoint union $A\sqcup B$ and whose order relation is the union of the order relation of~$A$, the order relation of~$B$, and the requirement that $a\leq b$ for all $a\in A$ and $b\in B$.
If $A_1, \ldots, A_r$ are finite posets, then $A_1*A_2* \ldots *A_r$ is defined inductively.\par

A {\em pole poset} is a poset of the form $A_1*A_2* \ldots *A_r$ where each $A_i$ either has cardinality one, or consists of exactly two incomparable elements.
If $a\in A_i$, then clearly $a$ has height~$i-1$ (with the usual convention that the minimal elements, that is, those in~$A_1$, have height~0).
The definition implies that there are two types of elements in a pole poset~$P$~:
\begin{enumerate}
\item If $A_i=\{a\}$ has cardinality one, then $a$ is comparable to every element of~$P$.
\item If $A_i=\{a,b\}$ has cardinality two, then $a$ is comparable to every element of $P-\{b\}$.
In that case, $b$ will be called the {\em twin} of~$a$ and written $\breve a$. In particular, $\breve{\breve a}=a$.
\end{enumerate}
Notice that a totally ordered poset is a pole poset (with no twins).
We write $P_1$ for the set of elements of the first type (the `totally ordered' part of~$P$)
and $P_2$ for the set of elements of the second type (the `twin' part of~$P$).\par

A {\em pole lattice} is a lattice whose underlying poset is a pole poset.
Whenever $a$ and $\breve a$ are incomparable elements of height~$i$ in a pole lattice~$P$, then they are both join-irreducible and meet-irreducible.
In this case, there is a single element of height $i-1$,
namely $r(a)=a\meet \breve a$, and a single element of height $i+1$, namely $s(a)=a\vee \breve a$.
Clearly $r(a)=r(\breve a)$ and $s(a)=s(\breve a)$.
Also, $\widehat0$ is the unique element of height~0 and $\widehat1$ is the unique element of maximal height.

\result{Lemma} \label{irred-pole}
Let $P$ be a pole lattice
and let $X=\{\,\widehat0\,\}\cup\{s(a) \mid a\in P_2\}$.
Then $P-X$ is the set of irreducible elements of~$P$.
\fresult

\pf This is easy and is left to the reader.
\endpf

We want to show that pole posets can be characterized by a condition which will be useful later in Section~\ref{Section-characterization}.
Recall that a relation~$R$ on a set~$X$ is a subset $R\subseteq X\times X$ and that the opposite relation $R\op$ is defined by~:
$$(x,y)\in R\op \iff (y,x)\in R \mpoint$$
Moreover, the product of two relations $S$ and $T$ is the relation defined by
$$ST:=\{ \, (z,x)\in X\times X \,\mid\, \exists\;y\in X \;\text{ such that } \; (z,y)\in S \,\text{ and }\, (y,x)\in T \,\} \mpoint$$
Let $\Sigma_X$ be the symmetric group of all permutations of~$X$.
Associated with a permutation $\sigma\in\Sigma_X$, there is the relation
$$\Delta_\sigma:=\{(\sigma(x),x)\in X\times X\mid x\in X\} \mpoint$$
In particular, we write $\Delta_X:=\Delta_{\Id}$ for the identity morphism of the object~$X$.
The map $\sigma\mapsto \Delta_\sigma$ is a monoid homomorphism and $\Delta_\sigma$ is invertible for every $\sigma\in\Sigma_X$.
The symmetric group $\Sigma_X$ acts on relations by conjugation~:
we write $R^\sigma =\Delta_{\sigma^{-1}} R \Delta_\sigma$ and $\ls\sigma R=\Delta_\sigma R\Delta_{\sigma^{-1}}$.\par

\result{Proposition} \label{characterization-pole}
Let $P$ be a finite poset and let $R\subseteq P\times P$ be its order relation (i.e. $(x,y)\in R \iff x\leq y$).
Let $\overline R=(P\times P)-R$.
The following are equivalent~:
\begin{enumerate}
\item $P$ is pole poset.
\item There exists a permutation $\tau$ of~$P$ such that
$$\forall \; x,y\in P \,,\; \text{ if } \, x\not\leq y \,, \; \text{ then } y\leq \tau(x) \mpoint$$
\item There exists a permutation $\tau$ of~$P$ such that $\overline R\op \Delta_{\tau^{-1}}\subseteq R$.
\end{enumerate}
Moreover, if (b) holds, then $\tau$ can be chosen to be an automorphism of the poset~$P$ and, in that case,
it is unique and it satisfies $\tau(a)=\breve a$ for all twins $a\in P_2$.
\fresult

\pf First note that the equivalence of (b) and~(c) follows immediately from the definitions, because
$$x\not\leq y \iff (y,x)\in \overline R\op$$
while we always have $(x,\tau(x))\in \Delta_{\tau^{-1}}$.\mpn

Suppose that (a) holds and define $\tau$ to be the permutation that preserves heights and satisfies $\tau(a)=\breve a$ for all twins $a\in P_2$.
Let $x,y\in P$ such that $x\not\leq y$.
If $x$ is the unique element of its height, then $x$ is comparable to all elements of~$P$ and $\tau(x)=x$.
It follows that $y< x=\tau(x)$.
If $x$ and $\breve x$ are distinct elements of the same height, i.e. twins, then $x$ is comparable to every element of $P-\{\breve x\}$.
Therefore, if $y\neq \breve x$, then $y<x$, hence also $y<\breve x=\tau(x)$, while if $y=\breve x$, then $y=\tau(x)$.
This proves that we get $y\leq \tau(x)$ in all cases, hence (b) holds.\mpn

We assume now that (b) holds and we want to prove~(a).
We proceed by induction on the size of~$P$, starting from the obvious case when $|P|=1$.
Suppose first that $P$ has at least two distinct maximal elements $w$ and~$z$.
Since $w\not\leq z$, we have $z\leq \tau(w)$ by~(b), hence $z=\tau(w)$ by maximality of~$z$. Similarly, $w=\tau(z)$.
Now if $x\not\leq w$, then $w\leq \tau(x)$ by~(b), hence $w=\tau(x)$, so that $x=\tau^{-1}(w)=z$.
In other words, if $x\neq w$ and $x\neq z$, then $x<w$. Similarly, if $x\neq w$ and $x\neq z$, then $x<z$.
Therefore $w$ and $z$ are the unique maximal elements of~$P$ and $P=Q*\{w,z\}$, where $Q=P-\{w,z\}$.\par

If $x,y\in Q$ and $x\not\leq y$, then $y\leq \tau(x)$ by~(b). But the permutation $\tau$ exchanges $w$ and~$z$, so it restricts to a permutation of~$Q$.
Therefore (b) holds for the poset~$Q$ and, by induction, $Q$ is a pole poset. It follows that $P$ is a pole poset, as required.\par

Suppose now that $P$ has a single maximal element~$w$. If $x\neq w$, then $w\not\leq x$ by maximality, hence $x\leq\tau(w)$ by~(b).
If $\tau(w)=w$, then $\tau$ restricts to a permutation of~$Q=P-\{w\}$ and again we are done by induction.\par

So we assume now that our single maximal element~$w$ satisfies $\tau(w)\neq w$.
The condition $x\leq\tau(w)$ obtained above means that $\tau(w)$ is the unique maximal element of~$P-\{w\}$.
Assume by induction that $w>\tau(w)>\ldots>\tau^i(w)$ and that $\tau^j(w)$ is the unique maximal element of~$P-\{w,\ldots,\tau^{j-1}(w)\}$, for every $j=1,\ldots,i$.
Then if $x\neq w,\tau(w),\ldots, \tau^i(w)$, we have $\tau^i(w)\not\leq x$, hence $x\leq \tau^{i+1}(w)$ by~(b).
But $\tau^{i+1}(w)\neq \tau(w),\ldots, \tau^i(w)$, otherwise $\tau^i(w) \in \{w,\tau(w),\ldots,\tau^{i-1}(w) \}$ which is impossible by our induction assumption.
Therefore, either $\tau^{i+1}(w)$ is the unique maximal element of~$P-\{w,\ldots,\tau^i(w)\}$ and we continue our induction argument, or $\tau^{i+1}(w)=w$.\par

Our induction argument must stop and we let $r\geq2$ be the smallest integer such that $\tau^r(w)=w$.
Then $w>\tau(w)>\ldots>\tau^{r-1}(w)$ and $\tau^j(w)$ is the unique maximal element of~$P-\{w,\ldots,\tau^{j-1}(w)\}$, for every $j=1,\ldots,r-1$.
Moreover, setting $Q=P-\{w,\ldots,\tau^{r-1}(w)\}$, we obtain
$$P=Q*\{\tau^{r-1}(w)\}*\ldots*\{\tau(w)\}*\{w\}$$
and $Q$ must be invariant under~$\tau$. By our main induction procedure, $Q$ is a pole poset.
It follows that $P$ is a pole poset. This proves~(a) and we are done.\mpn

In order to prove our additional statement, we continue the analysis of the permutation~$\tau$, as above.
In the case when $P$ has two maximal elements $w$ and $z$, then we have seen that $\tau(w)=z$.
Moreover $\tau$ restricts to a permutation of~$Q=P-\{w,z\}$.
By induction, $\tau_{\mid Q}$ can be replaced uniquely by an automorphism $\alpha$ of the pole poset~$Q$
such that $\alpha$ exchanges all the twins of~$Q$.
Extending $\alpha$ to~$P$ by requiring that $\alpha$ exchanges $w$ and~$z$ (as it must, as we have seen for~$\tau$),
we obtain an automorphism of~$P$ having the additional required properties.

In the case when $P$ has a single maximal element $w$, then we have seen that $\tau$ permutes cyclically the subset
$S=\{w,\tau(w),\ldots,\tau^{r-1}(w) \}$ for some $r\geq1$, and it restricts to a permutation of~$Q=P-S$.
By induction, $\tau_{\mid Q}$ can be replaced uniquely by an automorphism $\alpha$ of the pole poset~$Q$ such that $\alpha$ exchanges all the twins.
Extending $\alpha$ by the identity on~$S$, we obtain an automorphism of~$P$ having the additional required properties.
\endpf

A {\em join-morphism} from a lattice $T$ to a lattice~$T'$ is a map $f:T\to T'$ which commutes with joins, i.e. such that
$$f\big(\bigvee_{a\in A} a\big)=\bigvee_{a\in A}f(a)\mvirg$$
for any subset $A$ of~$T$.
Similarly, a {\em meet-morphism} is a map which commutes with meets.
It is easy to see that a join-morphism is order-preserving, by considering the join $t_1\vee t_2$ in the case where $t_1\leq_T t_2$ in the lattice~$T$.
Moreover, the case $A=\emptyset$ shows that a join-morphism maps $\widehat0\in T$ to~$\widehat0\in T'$.
The following result is well-known.

\result{Lemma} \label{distributive-extend}
Let $P$ and $T$ be finite lattices.
Suppose that $P$ is distributive and let $E=\Irr(P)$.
Then any order-preserving map $\varphi:E\to T$ extends uniquely to a join-morphism $\widetilde\varphi:P\to T$.
\fresult

\pf
For any $p\in P$, we can write uniquely
$$p=\bigvee_{\substack{e\in E \\ e\leq p}}e$$
and then define the extension of~$\varphi$ by
$$\widetilde\varphi(p)=\bigvee_{\substack{e\in E \\ e\leq p}} \varphi(e) \mpoint$$
To check that $\widetilde\varphi$ is a join-morphism, we use the fact that, for any $e\in E$ and $p,p'\in P$, we have
$$e\leq p\vee p' \iff e\leq p \;\text{ or } \; e\leq p' \mpoint$$
This is because, if $e\leq p\vee p'$, then, $e=e\wedge (p\vee p')=(e\wedge p)\vee(e\wedge p')$ by distributivity,
hence by irreducibility, either $e=e\wedge p$, i.e., $e\leq p$, or $e=e\wedge p'$, i.e., $e\leq p'$.
\endpf

\result{Notation} \label{notation-CL}
\begin{enumerate}
\item We let $\CL$ be the category whose objects are the finite lattices and
where, for any finite lattices $P$ and $T$, $\Hom_\CL(P,T)$ is the set of all join-morphisms from $P$ to~$T$.
\item We denote by $\Inj_\CL(P,T)$ the set of all injective join-morphisms $P\to T$.
\item We denote by $\Sur_\CL(T,P)$ the set of all surjective join-morphisms $T\to P$.
\end{enumerate}
\fresult

Recall from Section~8 of~\cite{BT3} that, for any join-morphism $f:T\to P$,
there is an opposite morphism $f\op:P\op \to T\op$ defined by
$$f\op(p)=\bigvee_{f(t)\leq p} t \,.$$

\result{Lemma} \label{op-morphisms}
Let $P$ and $T$ be finite lattices and let $f:T\to P$ be a join-morphism.
\begin{enumerate}
\item $f\op:P\op \to T\op$ is a join-morphism. In other words, for any subset $A$ of~$P$,
$$f\op(\bigwedge_{a\in A}a)=\bigwedge_{A\in A}f\op(a)$$
(because the meet~$\meet$ is the join in the opposite lattice).
\item If $g:P\to Q$ is a join-morphism, then $(gf)\op=f\op g\op$.
\item $(f\op)\op=f$.
\item If $f$ is surjective, then $ff\op=\Id_P$. In particular, $f\op$ is injective and, for any $p\in P$,
$$f\op(p)=\bigvee_{f(t)=p} t =\sup\big(f^{-1}(p)\big) \mpoint$$
\item If $f$ is injective, then $f\op f=\Id_T$. In particular, $f\op$ is surjective.
\item Passing to the opposite induces bijections $\Inj_\CL(P,T)\to \Sur_\CL(T\op,P\op)$
and $\Sur_\CL(T,P)\to \Inj_\CL(P\op,T\op)$.
\end{enumerate}
\fresult

\pf 
(a), (b) and (c) are proved in Section~8 of~\cite{BT3}.\mpn

(d) Let $p\in P$. The equality $\displaystyle\bigvee_{f(t)\leq p} t=\bigvee_{f(t)=p} t$ follows from the fact that $f$ is surjective and order-preserving.
Moreover, it is clear that $\displaystyle\bigvee_{f(t)=p} t=\sup\big(f^{-1}(p)\big)$ because $f$ is a join-morphism.
Finally $ff\op=\Id_P$ because $f\big( \sup\big(f^{-1}(p)\big)\big)=p$.\mpn

(e) This follows from (b), (c), and (d) by passing to opposite morphisms.\mpn

(f) This follows from (d) and (e).
\endpf

For later use, we now prove a specific result in the case when $P$ is a pole lattice.

\result{Lemma} \label{bijection-bases}
Let $T$ be a finite lattice and let $P$ be a pole lattice. Then there is a bijection between $\Inj_\CL(P,T)$ and $\Sur_\CL(T,P)$.
\fresult

\pf
Associated with the pole lattice $P$, there is the set
$$E_2=\{a_1,\breve{a}_1, a_2,\breve{a}_2, \;\ldots\; , a_n,\breve{a}_n \}$$
consisting of all the twins $a_i,\breve{a}_i$, indexed in such a way that $a_1<a_2<\ldots<a_n$.
Here $n$ is a positive integer (which is zero whenever $P$ is totally ordered).
We define
$$w_i=a_i\wedge \breve a_i \,, \qquad v_i=a_i\vee \breve a_i \,, \qquad  (1\leq i\leq n) \mvirg$$
and we also set $v_0=\widehat 0$ and $w_{n+1}=\widehat 1$.
Just above the pair of twins $a_i,\breve{a}_i$, there is a totally ordered interval $[v_i, w_{i+1}]$.
Also, we have a totally ordered interval $[v_0,w_1]$ below the pair $a_1,\breve{a}_1$,
and a totally ordered interval $[v_n,w_{n+1}]$ above the pair $a_n,\breve{a}_n$.
Note that we may have $v_i=w_{i+1}$.\par

Let $\lambda\in\Inj_\CL(P,T)$.
We want to define an injective meet-morphism $\widetilde\lambda:P\to T$ associated with~$\lambda$.
First we set
$$\widetilde\lambda(a_i)=\lambda(a_i) \,, \qquad \widetilde\lambda(\breve a_i)=\lambda(\breve a_i) \qquad(1\leq i\leq n) \mpoint$$
Since $\lambda$ is a join-morphism, we have
$$\lambda(v_0)=\widehat 0 \,, \qquad \lambda(v_i)=\lambda(a_i)\vee \lambda(\breve a_i)
\qquad(1\leq i\leq n) \mpoint$$
Note also that
$$\lambda(w_i)\leq \lambda(a_i)\wedge \lambda(\breve a_i)
\qquad(1\leq i\leq n) \,, \qquad \lambda(w_{n+1})\leq \widehat 1 \mpoint$$
We have to define $\widetilde\lambda$ on each interval $[v_{i-1},w_i]$, $1\leq i\leq n+1$, and there are two cases for each~$i$.\par

If $1\leq i\leq n$, either $\lambda(w_i)= \lambda(a_i)\wedge \lambda(\breve a_i)$ or
$\lambda(w_i)< \lambda(a_i)\wedge \lambda(\breve a_i)$.
In the first case, we simply set
$$\widetilde\lambda(x)=\lambda(x) \,,\qquad \forall\; x\in[v_{i-1},w_i] \mvirg$$
while in the second, we set
$$\widetilde\lambda(x)=\lambda(s(x)) \qquad \forall\; x\in[v_{i-1},w_i[ \,, \qquad \text{and } \qquad
\widetilde\lambda(w_i)=\lambda(a_i)\wedge \lambda(\breve a_i) \mvirg$$
where $s$ denotes the shift upwards in the totally ordered interval $[v_{i-1},w_i]$,
that is, $s(x)=\inf\{y\mid x<y\}$.\par

Similarly, if $i=n+1$, either $\lambda(w_{n+1})= \widehat 1$ or $\lambda(w_{n+1})< \widehat 1$.
In the first case, we simply set
$$\widetilde\lambda(x)=\lambda(x) \,,\qquad \forall\; x\in[v_n,w_{n+1}] \mvirg$$
while in the second, we set
$$\widetilde\lambda(x)=\lambda(s(x)) \qquad \forall\; x\in[v_n,w_{n+1}[ \,, \qquad \text{and } \qquad
\widetilde\lambda(w_{n+1})=\widehat 1 \mpoint$$
It is easy to see that $\widetilde\lambda$ is order-preserving and injective, and moreover
$$\widetilde\lambda(a_i\wedge \breve a_i )=\widetilde\lambda(w_i)=\lambda(a_i)\wedge \lambda(\breve a_i)
=  \widetilde\lambda(a_i)\wedge \widetilde\lambda(\breve a_i) \mpoint$$
In view of the structure of pole lattices, this means that $\widetilde\lambda:P\to T$ is a meet-morphism,
or in other words a join-morphism $\widetilde\lambda:P\op\to T\op$.
Therefore $\widetilde\lambda\in \Inj_\CL(P\op,T\op)$ and
this defines a map
$$\Omega_{P,T}: \Inj_\CL(P,T) \longrightarrow \Inj_\CL(P\op,T\op) \;,\qquad \lambda\mapsto\widetilde\lambda \mpoint$$

In the other direction, we proceed as follows.
The same construction, applied to $P\op$ and $T\op$, defines a map
$$\Omega_{P\op,T\op}: \Inj_\CL(P\op,T\op) \longrightarrow \Inj_\CL(P,T)$$
and it is elementary to check that $\Omega_{P\op,T\op}$ maps $\widetilde\lambda$ to~$\lambda$,
because the shift upwards $x\mapsto s(x)$ in the opposite $[v,w]\op$ of a totally ordered interval 
corresponds to the shift downwards $x\mapsto r(x)$ in the original interval $[v,w]$.
In other words the composite $\Omega_{P\op,T\op} \circ\Omega_{P,T}$ is the identity.
Similarly, $\Omega_{P,T} \circ \Omega_{P\op,T\op}$ is the identity and it follows that $\Omega_{P\op,T\op}$ is a bijection.\par

Now it suffices to compose with the bijection $\Inj_\CL(P\op,T\op) \to\Sur_\CL(T,P)$ of Lemma~\ref{op-morphisms}
to obtain a bijection between $\Inj_\CL(P,T)$ and $\Sur_\CL(T,P)$.
\endpf

%%%%%%%% Section

\section{Correspondence functors} \label{Section-correspondence-functors}

\bigskip
\noindent
We recall the basic facts we need about correspondence functors and we refer to Sections~2--4 of~\cite{BT2} and Section~2 of~\cite{BT3} for more details.
We denote by $\CC$ the category of finite sets and correspondences.
Its objects are the finite sets and the set $\CC(Y,X)$ of morphisms from $X$ to~$Y$ (using a reverse notation which is convenient for left actions)
is the set of all correspondences from $X$ to~$Y$, namely all subsets of $Y\times X$.
A correspondence from $X$ to~$X$ is also called a {\em relation} on~$X$.
Given two correspondences $R\subseteq Z\times Y$ and $S\subseteq Y\times X$,
their composition $RS$ is defined by
$$RS:=\{ \, (z,x)\in Z\times X \,\mid\, \exists\;y\in Y \;\text{ such that } \; (z,y)\in R \,\text{ and }\, (y,x)\in S \,\} \mvirg$$
and this generalizes the product of relations, defined in Section~\ref{Section-pole}.\par

For any commutative ring~$k$, we let $k\CC$ be the $k$-linearization of~$\CC$.
The objects are again the finite sets and $k\CC(Y,X)$ is the free $k$-module with basis $\CC(Y,X)$.
A {\em correspondence functor} is a $k$-linear functor from $k\CC$ to $k\text{-\!}\Mod$.
We let $\CF_k$ be the category of all correspondence functors (for some fixed commutative ring~$k$).
We define a {\em minimal set} for a correspondence functor~$F$ to be a finite set~$X$ of minimal cardinality such that $F(X)\neq \{0\}$.
For a nonzero functor, such a minimal set always exists and is unique up to bijection.\par

The first instances of correspondence functors are the representable functors $k\CC(-,E)$, where $E$ is a finite set, and more generally the functors
$$L_{E,W}:=k\CC(-,E)\otimes_{k\CC(E,E)}W$$
where $W$ is a left $k\CC(E,E)$-module.
Actually, the functor $W\mapsto L_{E,W}$ is left adjoint of the evaluation functor
$$\CF_k\longrightarrow k\CC(E,E)\,\text{-}\Mod \;, \qquad F\mapsto F(E) \mpoint$$
The correspondence functor $L_{E,W}$ has a subfunctor $J_{E,W}$ defined on any finite set~$X$ by
$$J_{E,W}(X):=\Big\{\sum_i\phi_i\otimes w_i \in L_{E,W}(X) \mid
\forall\psi\in k\CC(E,X), \sum_i (\psi\phi_i)\cdot w_i=0 \Big\} \mpoint$$
We shall work with the functor $L_{E,W}/J_{E,W}$ for some specific choices of $k\CC(E,E)$-modules~$W$.\par

Recall from~\cite{BT1} that, for a suitable two-sided ideal~$I$, there is a quotient algebra $\CP_E=k\CC(E,E)/I$,
called the {\em algebra of permuted orders} because it has a $k$-basis consisting of all relations on~$E$ of the form
$\Delta_\sigma R$, where $\sigma$ runs through the symmetric group $\Sigma_E$
and $R$ is an order relation on~$E$.
The product of two order relations $R$ and $S$ in~$\CP_E$ is the transitive closure of $R\cup S$ if this closure is an order relation, and zero otherwise.
This product, together with the conjugation action of permutations on relations, describes completely the algebra structure of~$\CP_E$.\par

Among the $k\CC(E,E)$-modules, there is the {\em fundamental module} $\CP_E f_R$,
associated with a poset $(E,R)$, where $E$ is a finite set and $R$ denotes the order relation on~$E$ which defines the poset structure.
Here $f_R$ is a suitable idempotent in~$\CP_E$, depending on~$R$, and $\CP_E f_R$ is the left ideal generated by~$f_R$.
The main thing we need to know about the fundamental module $\CP_E f_R$ is its structure as a bimodule.
This is described in the next result, which combines Corollary 7.3 and Proposition~8.5 of~\cite{BT1}.

\result{Proposition} \label{fundamental-module}
Let $E$ be a finite set and $R$ an order relation on~$E$.
\begin{enumerate}
\item The fundamental module $\CP_E f_R$ is a $(k\CC(E,E),k\Aut(E,R))$-bimodule and the right action of $k\Aut(E,R)$ is free.
\item $\CP_E f_R$ is a free $k$-module with a $k$-basis consisting of the elements $\Delta_\sigma f_R$,
where $\sigma$ runs through the group $\Sigma_E$ of all permutations of $E$.
\item The action of the algebra of relations $k\CC(E,E)$ on the module $\CP_E f_R$
is given as follows. For any relation $Q\in \CC(E,E)$,
$$Q\cdot\Delta_\sigma f_R=\left\{\begin{array}{ll}
\Delta_{\tau\sigma}f_R&\hbox{if}\;\;\exists\,\tau\in\Sigma_E\;\hbox{such that}\;
\Delta_E\subseteq \Delta_{\tau^{-1}}Q\subseteq {\ls\sigma R},\\
0&\hbox{otherwise}\mpoint\end{array}\right.$$
Moreover, $\tau$ is unique in the first case.
\end{enumerate}
\fresult

Using the bimodule structure on~$\CP_E f_R$, we define
$$T_{R,V}:=\CP_E f_R\otimes_{k\Aut(E,R)} V \mvirg$$
where $V$ is any $k\Aut(E,R)$-module.
Then $T_{R,V}$ is a left $k\CC(E,E)$-module for the action induced from the action of $k\CC(E,E)$ on~$\CP_Ef_R$
described in Proposition~\ref{fundamental-module} above.
The main thing we need to know about $T_{R,V}$ is the following result, which is part of Theorem~8.1 in~\cite{BT1}.

\result{Proposition} \label{simple-modules} Assume that $k$ is a field.
If $V$ is a simple $k\Aut(E,R)$-module, then $T_{R,V}$ is a simple $k\CC(E,E)$-module.\fresult

Associated with the above $k\CC(E,E)$-modules, we can now define, as in~\cite {BT2} and~\cite{BT3}, some specific correspondence functors.
Using the fundamental module $\CP_E f_R$, we define
$$\S_{E,R}:=L_{E,\CP_E f_R}/J_{E,\CP_E f_R}$$
and we call it the {\em fundamental functor} associated with the poset $(E,R)$.
Using the module $T_{R,V}$, we define
$$S_{E,R,V}:=L_{E,T_{R,V}}/J_{E,T_{R,V}} \mpoint$$
Note that $\S_{E,\ls\sigma R}\cong \S_{E,R}$ and $S_{E,\ls\sigma R,\ls\sigma V}\cong S_{E,R,V}$, for any permutation $\sigma\in\Sigma_E$.
Our next result is Proposition~2.6 in~\cite{BT3}.

\result{Proposition} \label{SER-SERV}
\begin{enumerate}
\item The set $E$ is a minimal set for $\S_{E,R}$ and $\S_{E,R}(E)\cong \CP_E f_R$ as left $k\CC(E,E)$-modules.
\item The set $E$ is a minimal set for $S_{E,R,V}$ and $S_{E,R,V}(E)\cong T_{R,V}$ as left $k\CC(E,E)$-modules.
\item If $k$ is a field and $V$ is a simple $k\Aut(E,R)$-module, then $S_{E,R,V}$ is a simple correspondence functor.
\end{enumerate}
\fresult

It is proved in Theorem~4.7 of~\cite{BT2} that, when $k$ is a field, any simple functor has the form $S_{E,R,V}$ for some triple $(E,R,V)$ and that
the set of isomorphism classes of simple correspondence functors is parametrized by
the set of isomorphism classes of triples $(E,R,V)$ where $E$ is a finite set, $R$ is an order relation on~$E$, and $V$ is a simple $k\Aut(E,R)$-module.\par

We note that the fundamental functor $\S_{E,R}$ is a precursor of~$S_{E,R,V}$, in the sense of the following result.

\result{Proposition} \label{precursor}
Suppose that $V$ is a simple $k\Aut(E,R)$-module, hence in particular generated by a single element~$v$.
\begin{enumerate}
\item Consider the surjective morphism of correspondence functors
$$\Phi: L_{E,\CP_Ef_R} \longrightarrow L_{E,T_{R,V}}$$
induced by the surjective homomorphism of $\CP_E$-modules
$$\Phi_E:\CP_Ef_R \longrightarrow \CP_Ef_R \otimes_{k\Aut(E,R)} V =T_{R,V} \,, \qquad a\mapsto a\otimes v \mpoint$$
Then $\Phi$ induces a surjective morphism of correspondence functors
$$\S_{E,R} \longrightarrow S_{E,R,V} \mpoint$$
\item $\Phi$ induces an isomorphsim
$$\S_{E,R}\otimes_{k\Aut(E,R)} V \cong S_{E,R,V} \mpoint$$
\end{enumerate}
\fresult

\pf
(a) is Lemma~2.7 in~\cite{BT3}, while (b), which is far from being obvious, is Theorem~6.10 in~\cite{BT4}.
\endpf

In short, it it is possible to recover $S_{E,R,V}$ from~$\S_{E,R}$ by simply tensoring with~$V$.
Consequently, the fundamental functors play a crucial role throughout our work.\mpn

Another important construction of correspondence functors is obtained from finite lattices (see~\cite{BT3} for details).

\result{Definition} \label{def-F}
If $T$ is a finite lattice and $X$ is a finite set, define $F_T(X)=k T^X$, the free $k$-module on the set~$T^X$ of all functions $X\to T$.
Given $\varphi: X\to T$ and a correspondence $S\in\CC(Y,X)$, then $S\varphi: Y\to T$ is defined by the formula
$$(S\varphi)(y) =\bigvee_{(y,x)\in S} \varphi(x) \mpoint$$
Then $F_T$ becomes in this way a correspondence functor.
\fresult

We want to recall two main properties of this construction but we first need some notation.
Let $\CL$ be the category of finite lattices and join-morphisms, as in Notation~\ref{notation-CL}.
The $k$-linearization $k\CL$ of~$\CL$ has the same objects
and $\Hom_{k\CL}(T,T')$ is the free $k$-module $k\Hom_\CL(T,T')$ with basis $\Hom_\CL(T,T')$.
The composition of morphisms in $k\CL$ is the $k$-bilinear extension of the composition in~$\CL$.
The following results appear in Theorems~4.8 and 4.12 of~\cite{BT3}.

\result{Theorem} \label{fully-faithful}
\begin{enumerate}
\item The assignment $T\mapsto F_T$ extends to a $k$-linear functor $F_?:k\CL \to \CF_k$. Moreover, $F_?$ is fully faithful.
\item If $T$ is a finite lattice, then $F_T$ is projective in~$\CF_k$ if and only if $T$ is distributive.
In particular, if $P$ is a pole lattice, then $F_P$ is projective.
\end{enumerate}
\fresult

Our next lemma gives another realization of the functor~$F_T$ in a special case.
Let $E$ be a finite set and $R$ an order relation on~$E$ (i.e. $(E,R)$ is a finite poset).
As  in~\cite{BT3}, let $\Idown(E,R)$ be the lattice of all subsets of~$E$ closed under taking smaller elements with respect to~$R$.
Then $(E,R)$ is isomorphic to the poset of irreducible elements of~$\Idown(E,R)$ via the map $e\mapsto E_{\leq e}=\{x\in E\mid x\leq e\}$.
Notice that $r(E_{\leq e})=E_{< e}$ in the lattice $\Idown(E,R)$.

\result{Lemma} \label{iso-complement}
Let $(E,R)$ be a finite poset and let $T=\Idown(E,R)$.
For any finite set~$X$, define a map
$$\rho_X:F_{T\op}(X) \longrightarrow k\CC(X,E)R \,,\qquad
\rho_X(\varphi)= \{(x,e) \mid e\notin \varphi(x)\} \subseteq X\times E \mvirg$$
where $\varphi:X\to T\op$ is any basis element in~$F_{T\op}(X)$.
Then this induces an isomorphism of correspondence functors
$\rho:F_{T\op}\longrightarrow k\CC(-,E)R$.
\fresult

\pf
The result can be obtained by combining Proposition~4.5 and Remark~8.7 in~\cite{BT3}, using the isomorphism, via complementation,
$\Idown(E,R\op) \cong \Idown(E,R)\op$. We provide instead a direct proof.\par

Since $\rho_X(\varphi)$ is a subset of~$X\times E$, it is an element of~$\CC(X,E)$.
It is right invariant by~$R$
because if $(x,e)\in\rho_X(\varphi)$, i.e.~$e\notin\varphi(x)$, and if $(e,f)\in R$,
then $(x,f)\in \rho_X(\varphi)$ because $f\notin \varphi(x)$
(otherwise we would have $e\in\varphi(x)$ since $\varphi(x)$ is closed under taking smaller elements).
Hence $\rho_X(\varphi)=\rho_X(\varphi)R\in \CC(X,E)R$.
It is elementary to check that $\rho$ is a morphism of functors.
Moreover, it is an isomorphism because there is an inverse morphism mapping $S\in \CC(X,E)R$ to the function $\varphi_S:X\to T\op$ defined by
$$\varphi_S(x)=\{e\in E \,\mid\, (x,e)\notin S \} \mpoint$$
The fact that $S$ is right invariant by~$R$ implies that $\varphi_S(x)$ is closed under taking smaller elements.
Details are left to the reader.
\endpf

There is a direct connection between the functors associated with lattices and the fundamental functors.
This is Theorem~6.5 in~\cite{BT3}.

\result{Theorem} \label{Theta}
Let $(E,R)$ be a finite poset.
There is a unique surjective morphism
$$\Theta:F_{\Idown(E,R\op)} \longrightarrow \S_{E,R}$$
mapping the inclusion map $j: E \to \Idown(E,R\op)$ to $f_R\in \S_{E,R}(E)\cong \CP_Ef_R$.
\fresult

We now recall another main result from~\cite{BT3}, which will be used in Section~\ref{Section-characterization} and Section~\ref{Section-fundamental}.
Let $T=\Idown(E,R)$ and, as in Section~9 of~\cite{BT3}, consider the element
\begin{equation} \label{def-gamma}
\gamma_T:=\sum_{A\subseteq E} (-1)^{|A|} \eta_A^0 \in F_{T\op}(E) \mpoint
\end{equation}
Here $\eta_A^0:E\to T\op$ is the map defined by
$$\eta_A^0(e)= \left\{\begin{array}{ll}
r(E_{\leq e})=E_{<e} &\hbox{if}\;\; e\in A \mvirg \\
E_{\leq e} &\hbox{if}\;\; e\notin A \mvirg
\end{array}\right.$$
with values in the lattice~$T$, but viewed as elements of~$T\op$.

\result{Theorem} \label{SERgamma}
Let $(E,R)$ be a finite poset and let $T=\Idown(E,R)$.
The subfunctor of~$F_{T\op}$ generated by $\gamma_T$ is isomorphic to the fundamental functor~$\S_{E,R}$.
Moreover, the isomorphism
$${<}\gamma_T{>}(E) \longrightarrow \S_{E,R}(E) \cong \CP_E f_R$$
maps $\gamma_T$ to~$f_R$.
\fresult

\pf
The first statement is Theorem~9.5 in~\cite{BT3}. The second statement can be traced in the proof of that theorem.
More precisely, if  $j:E \to \Idown(E,R\op)$ denotes the inclusion map,
it is shown that $\gamma_T\in F_{T\op}(E)$ is the image of $j\in F_{\Idown(E,R\op)}(E)$ under a morphism
$$\xi: F_{\Idown(E,R\op)} \longrightarrow F_{T\op} \mpoint$$
On the other hand, by Theorem~\ref{Theta} above, there is a surjective morphism
$$\Theta:F_{\Idown(E,R\op)} \longrightarrow \S_{E,R}$$
mapping the inclusion map $j$ to $f_R\in \S_{E,R}(E)\cong \CP_Ef_R$.
Both morphisms $\xi$ and $\Theta$ are proved to have the same kernel and this induces the required isomorphism ${<}\gamma_T{>} \cong \S_{E,R}$.
It follows that this isomorphism maps $\gamma_T$ to~$f_R$.
\endpf

%%%%%%%% Section

\section{Characterization of simple projective functors} \label{Section-characterization}

\bigskip
\noindent
Throughout this section, assume that the base ring $k$ is a field and let $(E,R)$ be a finite poset.
Our aim is to characterize the triples $(E,R,V)$ such that the simple correspondence functor $S_{E,R,V}$ is projective.\par

Since $S_{E,R,V}$ is isomorphic to a quotient of the fundamental functor~$\S_{E,R}$ (see Proposition~\ref{precursor}), we shall actually work with the latter.
We have $\S_{E,R}\cong {<}\gamma_T{>}$ by Theorem~\ref{SERgamma}, where $T=\Idown(E,R)$ and $\gamma_T$ is defined by~(\ref{def-gamma}).
We let
$$\zeta: \,{<}\gamma_T{>} \,\longrightarrow F_{T\op}$$
be the inclusion morphism.
We also let
$$\rho:F_{T\op}\longrightarrow k\CC(-,E)R$$
be the isomorphism of correspondence functors described in Lemma~\ref{iso-complement} and we define
$$\delta:=\rho\zeta(\gamma_T)=\rho(\gamma_T) \in k\CC(E,E)R \mpoint$$
In view of the isomorphism~$\rho$, the subfunctor ${<}\delta{>}$ of~$k\CC(-,E)R$ generated by~$\delta$ is isomorphic to~${<}\gamma_T{>}$, hence to~$\S_{E,R}$.
We shall work with~$\delta$ and we first need its precise description as a linear combination of relations.

\result{Lemma} \label{delta}
Let $\delta:=\rho(\gamma_T) \in k\CC(E,E)R$, where $\gamma_T$ is defined by~(\ref{def-gamma}).
\begin{enumerate}
\item $\rho(\eta_A^0)=\overline R\op \cup \Delta_A$, where $\overline R=(E\times E)-R$ and $\Delta_A =\{ (a,a)\mid a\in A \}\subseteq E\times E$.
\item $\delta=\displaystyle\sum_{A\subseteq E} (-1)^{|A|} (\overline R\op \cup \Delta_A)$.
\item $R(\overline R\op \cup \Delta_A)=\overline R\op \cup \Delta_A$.
\item $R\gamma_T=\gamma_T$ and $R\delta=\delta$.
\end{enumerate}
\fresult

\pf Throughout this proof, we write $x\leq y$ for $(x,y)\in R$.\mpn

(a) By Lemma~\ref{iso-complement}, we have
$$\rho(\eta_A^0)=\{ (f,e)\in E\times E \mid e\notin \eta_A^0(f) \}
= \left\{\begin{array}{ll}
\{ (f,e)\in E\times E \mid e\not< f \} &\hbox{if}\;\; e\in A \mvirg \\
\{ (f,e)\in E\times E \mid e\not\leq f \} &\hbox{if}\;\; e\notin A \mpoint
\end{array}\right.$$
But $\{ (f,e)\in E\times E \mid e\not\leq f \}=\overline R\op$.
If $e\in A$, we need to add to~$\overline R\op$ the element $(e,e)$, because $e\not< e$. Therefore
$\rho(\eta_A^0)=\overline R\op \cup \Delta_A$, as required.\mpn

(b) This follows from (a) and the fact that $\delta=\rho(\gamma_T)=\sum_{A\subseteq E} (-1)^{|A|} \rho(\eta_A^0)$.\mpn

(c) Since $\Delta_E\subseteq R$, we have an inclusion
$$\overline R\op \cup \Delta_A = \Delta_E(\overline R\op \cup \Delta_A) \subseteq R(\overline R\op \cup \Delta_A) \mpoint$$
In order to prove the reverse inclusion, we let $(a,c)\in R(\overline R\op \cup \Delta_A)$.
Then there exists $b\in E$ such that $a\leq b$ and $(b,c)\in\overline R\op \cup \Delta_A$.\par

If $(b,c)\in\overline R\op$, that is, $c\not\leq b$, then $c\not\leq a$, otherwise we would have $c\leq a \leq b$.
Therefore $(a,c)\in\overline R\op$.\par

If $(b,c)\in\Delta_A$, then $b=c\in A$ and there are two cases. If $a=b$, then $(a,c)=(a,a)\in\Delta_A$.
If $a\neq b$, then $a<b=c$, hence $c\not\leq a$, that is, $(a,c)\in \overline R\op$.\par

This completes the proof that $R(\overline R\op \cup \Delta_A)\subseteq \overline R\op \cup \Delta_A$, hence equality.\mpn

(d) It follows from (b) and (c) that $R\delta=\delta$, hence also $R\gamma_T=\gamma_T$
because $\rho$ is an isomorphism mapping $\gamma_T$ to~$\delta$.
The latter equality was also proved in Lemma~9.3 of~\cite{BT3}.
\endpf

We also need some technical computations involving~$\delta$.

\result{Lemma} \label{S-delta}
As above, consider $\delta=\displaystyle\sum_{A\subseteq E} (-1)^{|A|} (\overline R\op \cup \Delta_A)$.
Let $S\in\CC(E,E)R$ (that is, $S\subseteq E\times E$ and $S=SR$).
\begin{enumerate}
\item $S\delta\neq0$ if and only if
there exists a permutation $\sigma\in \Sigma_E$ such that $S=\Delta_\sigma R$.
\item If $S=RS$ and $S\delta\neq0$, then there exists an automorphism $\sigma\in \Aut(E,R)$ such that 
$S=\Delta_\sigma R$.
\item If $(E,R)$ is a pole poset and if $(\overline R\op \cup \Delta_A)\delta\neq0$, then
$$A=E_1 \qquad \text{ and } \qquad \overline R\op \cup \Delta_A=\Delta_\tau R \mvirg$$
where $\tau$ is the automorphism of $(E,R)$
satisfying $\tau(a)=\breve a$ for all $a\in E_2$ (the twin part of~$E$) and $\tau(a)=a$ for all $a\in E_1$ (the totally ordered part of~$E$).
\end{enumerate}
\fresult

\pf (a) The condition $S\delta\neq0$ is equivalent to $Sf_R\neq0$ by Theorem~\ref{SERgamma}.
By Proposition~\ref{fundamental-module}, we obtain
$$\begin{array}{rcl}
S\delta\neq0 \,\iff\, Sf_R\neq0 &\Longrightarrow& \exists\,\sigma\in\Sigma_E \;\text{ such that } \;\Delta_E \subseteq \Delta_\sigma^{-1}S\subseteq R \\
&\Longrightarrow& \exists\,\sigma\in\Sigma_E \;\text{ such that } \;R \subseteq \Delta_\sigma^{-1}SR\subseteq R^2 \\
&\Longrightarrow& \exists\,\sigma\in\Sigma_E \;\text{ such that } \;R \subseteq \Delta_\sigma^{-1}S\subseteq R \\
&\Longrightarrow& \exists\,\sigma\in\Sigma_E \;\text{ such that } S=\Delta_\sigma R \mvirg
\end{array}$$
using the equalities $S=SR$ and $R^2=R$ (by transitivity and reflexivity of~$R$).
Conversely, if $S=\Delta_\sigma R$, then, by Lemma~\ref{delta}, we obtain
$$S\delta=\Delta_\sigma R\delta=\Delta_\sigma \delta \neq0 \mvirg$$
because $\delta\neq0$ since it generates a nonzero subfunctor.\mpn

(b) We have $S=\Delta_\sigma R$ by~(a) and since $S=RS$, we obtain $R\Delta_\sigma R=\Delta_\sigma R$, or in other words
$R^\sigma R=R$, where $R^\sigma=\Delta_\sigma^{-1} R\Delta_\sigma$.
Since $\Delta_E\subseteq R$, we get $R^\sigma \subseteq R^\sigma R= R$, hence $R^\sigma = R$ because both relations have the same cardinality.
This means that $\Delta_\sigma$ commutes with~$R$, that is, $\sigma$ is an automorphism of the poset $(E,R)$.\mpn

(c) By (b) applied to $S=\overline R\op \cup \Delta_A$ (which satisfies $S=RS$ by Lemma~\ref{delta}), we have
$$\overline R\op \cup \Delta_A=\Delta_\sigma R$$
for some automorphism $\sigma\in \Aut(E,R)$.
Since $(E,R)$ is a pole poset, $\sigma$ is the identity on~$E_1$ and interchanges some of the twins $e,\breve e\in E_2$, so in particular $\sigma=\sigma^{-1}$.\par

If $e\in E_2$ and $\breve e$ is its twin, then $\breve e\not\leq e$, hence $(e,\breve e)\in\overline R\op\subseteq \Delta_\sigma R$.
Therefore $(e,\sigma(e))\in\Delta_\sigma$ and $(\sigma(e),\breve e)\in R$, that is $\sigma(e)\leq\breve e$.
This shows that $\sigma(e)$ cannot be equal to~$e$, i.e. $\sigma(e)=\breve e$.
Thus $\sigma$ interchanges all the twins, that is, it is equal to the automorphism~$\tau$ of the statement.\par

If $e\in E_1$, then $(e,e)\in \Delta_\sigma$ and $(e,e)\in R$, so $(e,e)\in \Delta_\sigma R$.
If conversely $(e,e)\in \Delta_\sigma R$, then $(e,\sigma(e))\in \Delta_\sigma$ and $(\sigma(e),e)\in R$, that is, $\sigma(e)\leq e$.
This cannot hold if $e\in E_2$, because $\sigma(e)=\breve e\not\leq e$, and therefore $e\in E_1$.
It follows that
$$e\in E_1 \iff (e,e)\in \Delta_\sigma R \iff (e,e) \in \overline R\op \cup \Delta_A
\iff (e,e) \in \Delta_A \iff e\in A \mvirg$$
the third equivalence using the fact that $(e,e)\notin \overline R\op$ because $e\leq e$.
This shows that $A=E_1$ and completes the proof.
\endpf

One of the key part of the proof of the main result is contained in the next lemma,
which will also be used again in Section~\ref{Section-fundamental}.

\result{Lemma} \label{proj-pole}
Suppose that $k$ is a field.
Let $\S_{E,R}$ be the fundamental functor associated with a finite poset $(E,R)$ and let $M$ be a direct summand of~$\S_{E,R}$.
If $M$ is projective, then $(E,R)$ is a pole poset.
\fresult

\pf
Since $\S_{E,R}\cong {<}\gamma_T{>}$ by Theorem~\ref{SERgamma},  we can view $M$ as a direct summand of~${<}\gamma_T{>}$
and we let $\omega:M \longrightarrow \,{<}\gamma_T{>}$ be the inclusion morphism.
As above, we let $\zeta: \,{<}\gamma_T{>} \,\longrightarrow F_{T\op}$ be the inclusion morphism and
$\rho:F_{T\op}\longrightarrow k\CC(-,E)R$
be the isomorphism of correspondence functors described in Lemma~\ref{iso-complement}.
Finally let
$$\alpha: M \,\longrightarrow k\CC(-,E)R$$
be the composite $\alpha=\rho\zeta \omega$.\par

Since $M$ is projective and the base ring $k$ is a field, $M$ is also injective, by Theorem~10.6 in~\cite{BT2}.
Therefore the injective morphism $\alpha$ splits, that is, there exists a surjective morphism
$$\sigma:k\CC(-,E)R \longrightarrow M$$
such that $\sigma\alpha=\Id$. Thus $\alpha\sigma$ is an idempotent endomorphism of~$k\CC(-,E)R$.
Since $R\in k\CC(E,E)$ is a generator of~$k\CC(-,E)R$, its image $c:=\sigma(R)\in M(E)$ is a generator of~$M$.
Now $\gamma_T$ generates  ${<}\gamma_T{>}$, so we can write $\omega(c)=v\gamma_T$ for some $v\in k\CC(E,E)$.
We know that $R\gamma_T=\gamma_T$ by Lemma~\ref{delta} and therefore $v\gamma_T=vR\gamma_T$.
Replacing $v$ by $vR$, we can assume that $v=vR$ and we do so. Thus $v\in k\CC(E,E)R$.
Note that $c\neq0$, hence $v\gamma_T\neq0$.\par

Now for any $u\in k\CC(X,E)$, we have
$$\alpha\sigma(uR)=u{\cdot} \alpha\sigma(R)=u{\cdot} \alpha(c)=u{\cdot} \rho\zeta \omega(c)
=u{\cdot} \rho\zeta(v\gamma_T)
=uv{\cdot} \rho\zeta(\gamma_T)=uv\delta \mvirg$$
where $\delta=\rho\zeta(\gamma_T)=\rho(\gamma_T)$ as in Lemma~\ref{delta}.
In particular, using the fact that $\delta=\delta R$ (because $\delta\in k\CC(E,E)R$), we obtain
$$\alpha\sigma(\delta)=\alpha\sigma(\delta R)=\delta v \delta \mpoint$$
Since $v\gamma_T$ is nonzero, so is its image $v\delta=\rho\zeta(v\gamma_T)$ under the injective morphism~$\rho\zeta$ and therefore
$$0\neq v\delta=\alpha\sigma(R)=(\alpha\sigma)^2(R)=\alpha\sigma(v\delta)=v{\cdot}\alpha\sigma(\delta)=v\delta v\delta \mvirg$$
from which it follows that $\delta v\delta\neq0$.\par

Summarizing, we have proved that, under the assumption that $M$ is projective,
the element $\delta=\rho(\gamma_T)\in k\CC(E,E)R$ satisfies~:
\begin{equation} \label{nonzero-condition}
\exists \, v\in k\CC(E,E)R \;\text{ with } \; \delta v\delta\neq0 \mpoint
\end{equation}
Our aim is to show that (\ref{nonzero-condition}) implies that $(E,R)$ is a pole poset.\par

The condition $\delta v\delta\neq0$ implies that there exists a relation~$S$
(in the expression of $v\in k\CC(E,E)R$ as a linear combination of relations) such that $\delta S\delta\neq0$.
In particular $S\delta\neq0$, hence $S=\Delta_\tau R$ for some $\tau\in\Sigma_E$, by Lemma~\ref{S-delta}.
In view of the expression of~$\delta$ obtained in Lemma~\ref{delta}, there exists a subset $A\subseteq E$ such that
$$(\overline R\op \cup \Delta_A) \Delta_\tau R \delta \neq0 \mpoint$$
Again, this implies that the relation $(\overline R\op \cup \Delta_A) \Delta_\tau R$ has the form
$$(\overline R\op \cup \Delta_A) \Delta_\tau R=\Delta_\sigma R$$
for some $\sigma\in \Sigma_E$, by Lemma~\ref{S-delta}.
Since the left hand side is invariant under left multiplication by~$R$ (by Lemma~\ref{delta}),
part~(b) of Lemma~\ref{S-delta} implies that $\Delta_\sigma$ commutes with~$R$ (i.e. $\sigma\in\Aut(E,R)$).
It follows that
$$(\overline R\op \cup \Delta_A) \Delta_{\tau\sigma^{-1}} R= R \mpoint$$
In particular, we deduce that
$$\overline R\op \Delta_\psi\subseteq R \,, \qquad\text{where } \;\psi:=\tau\sigma^{-1} \mpoint$$
By the characterization of Proposition~\ref{characterization-pole}, this implies that $(E,R)$ is a pole poset, as was to be shown.
\endpf

\result{Theorem} \label{characterization-proj}
Let $k$ be a field and let $S_{E,R,V}$ be the simple correspondence functor parametrized by a finite set~$E$, an order relation~$R$ on~$E$,
and a simple $k\Aut(E,R)$-module~$V$.
The following conditions are equivalent.
\begin{enumerate}
\item The simple correspondence functor $S_{E,R,V}$ is projective.
\item $(E,R)$ is a pole poset and $V$ is a projective $k\Aut(E,R)$-module.
\item Either $(E,R)$ is a totally ordered poset or $(E,R)$ is a pole poset and the characteristic of~$k$ is different from~2.
\end{enumerate}
\fresult

\pf (b) $\Longleftrightarrow$ (c). For a pole poset $(E,R)$, the group $\Aut(E,R)$ is a 2-group (elementary abelian), generated by all the possible transpositions of twins.
In case $(E,R)$ is totally ordered, this group is trivial and the unique simple $k$-module~$k$ is automatically projective.
In case $(E,R)$ is a pole poset but is not totally ordered, then $\Aut(E,R)$ is nontrivial and the characteristic of~$k$ comes into play.
If ${\rm char}(k)\neq2$, then all simple $k\Aut(E,R)$-module~$V$ are projective (by Maschke's theorem).
If ${\rm char}(k)=2$, then the unique simple $k\Aut(E,R)$-module is the trivial module, which is not projective (by the converse of Maschke's theorem).\mpn

(a) $\Rightarrow$ (b). Since $S_{E,R,V}$ is projective by assumption and isomorphic to a quotient of the fundamental functor~$\S_{E,R}$
by Proposition~\ref{precursor}, it is isomorphic to a direct summand of~$\S_{E,R}$.
Therefore Lemma~\ref{proj-pole} can be applied and it follows that $(E,R)$ is a pole poset.\par

We also have to prove that $V$ is a projective $k\Aut(E,R)$-module. Let
$$T=T_{R,V}=\CP_E f_R\otimes_{k\Aut(E,R)}V$$
be the simple $k\CC(E,E)$-module appearing in the definition $S_{E,R,V}:=L_{E,T}/J_{E,T}$.
By adjunction, there is an isomorphism
$$\End_{\CF_k}(L_{E,T}) \cong \Hom_{k\CC(E,E)}(T,L_{E,T}(E)) \cong \End_{k\CC(E,E)}(T)$$
and this is a skew field by Schur's lemma (it is actually the field~$k$).
This has no nontrivial idempotent and so $L_{E,T}$ is indecomposable.
But the surjective morphism
$$\pi: L_{E,T} \longrightarrow L_{E,T}/J_{E,T}=S_{E,R,V}$$
is split because $S_{E,R,V}$ is projective by assumption.
Therefore $\pi$ is an isomorphism, by indecomposability of~$L_{E,T}$, hence $L_{E,T}$ is projective.\par

Evaluating this projective functor at the finite set~$E$, we obtain a $k\CC(E,E)$-module
$$L_{E,T}(E)=T=\CP_E f_R\otimes_{k\Aut(E,R)}V$$
which must be projective, by Lemma~10.1 in~\cite{BT2}.
Now $\CP_E f_R\otimes_{k\Aut(E,R)}V$ is actually a module for the quotient algebra $\CP_E=k\CC(E,E)/I$
(see Section~\ref{Section-correspondence-functors}).
It follows that $\CP_E f_R\otimes_{k\Aut(E,R)}V$ is a projective module for the algebra~$\CP_E$,
because of the splitting of the composition of surjective homomorphisms
$$k\CC(E,E) \longrightarrow \CP_E \longrightarrow \CP_E f_R\otimes_{k\Aut(E,R)}V \mpoint$$

Finally, by Theorem~7.5 in~\cite{BT1}, there is an isomorphism of algebras
$$\CP_E \cong \prod_{R} M_{n_R}(k\Aut(E,R))$$
for some integers~$n_R$, where $R$ runs over all order relations on~$E$ up to isomorphism.
Thus there is a Morita equivalence
$$\CP_E\,\text{-}\Mod \;\cong\;  \prod_{R} \; k\Aut(E,R)\,\text{-}\Mod$$
and the bimodule inducing the equivalence is $\bigoplus_R \CP_E f_R$ (see Remark~7.6 in~\cite{BT1}).
Therefore $\CP_E f_R\otimes_{k\Aut(E,R)}V$ corresponds to the $k\Aut(E,R)$-module~$V$ under this equivalence.
Since projectivity is preserved by a Morita equivalence, $V$ is a projective $k\Aut(E,R)$-module, as required.\mpn

(b) $\Rightarrow$ (a). We assume that $(E,R)$ is a pole poset and, as before, we write $x\leq y$ for $(x,y)\in R$.
Our aim is to compute $\delta^2$ and show that it is an idempotent.
In view of the expression of~$\delta$ in Lemma~\ref{delta}, we have to consider terms of the form $(\overline R\op \cup \Delta_A)\delta$.
By Lemma~\ref{S-delta}, this can be nonzero only if $A=E_1$ and $\overline R\op \cup \Delta_A=\Delta_\tau R$,
where $\tau\in\Aut(E,R)$ is the automorphism exchanging all twins $e,\breve e\in E_2$
and fixing $E_1=E-E_2$ pointwise.\par

Thus $\overline R\op \cup \Delta_{E_1}$ is the only term which can come into play and
we now show that it is indeed equal to~$\Delta_\tau R$.
For any $a\in E_1$, we have $(a,\tau(a))\in \Delta_\tau$ and $\tau(a)=a\leq a$, hence $(a,a)\in\Delta_\tau R$.
Therefore $\Delta_{E_1}\subseteq \Delta_\tau R$.
Since $(E,R)$ is a pole poset, Proposition~\ref{characterization-pole} implies that $\overline R\op \subseteq R \Delta_\tau= \Delta_\tau R$,
using the fact that $\tau$ is an automorphism of~$(E,R)$.
So we obtain
$$\overline R\op \cup \Delta_{E_1} \subseteq \Delta_\tau R \mpoint$$
In order to prove the reverse inclusion, we let $(a,\tau(a))\in \Delta_\tau$ (since $\tau=\tau^{-1}$) and $(\tau(a),b)\in R$, i.e. $\tau(a)\leq b$.
If $a\in E_2$, then $\tau(a)=\breve a$, hence $\breve a \leq b$. Then $b\not\leq a$, that is, $(a,b)\in\overline R\op$.
If $a\in E_1$, then $\tau(a)=a$, hence $a\leq b$. If $a=b$, then $(a,b)\in \Delta_{E_1}$,
while if $a\neq b$, then $a<b$, hence $b\not\leq a$, that is, $(a,b)\in\overline R\op$.
This shows the required reverse inclusion and therefore
$$\overline R\op \cup \Delta_{E_1} = \Delta_\tau R \mvirg$$
as claimed.
In particular $(\overline R\op \cup \Delta_{E_1})\delta=\Delta_\tau R\delta=\Delta_\tau \delta$ by Lemma~\ref{delta}.
Therefore
$$\delta^2=\sum_{A\subseteq E} (-1)^{|A|} (\overline R\op \cup \Delta_A)\delta
=(-1)^{|E_1|}(\overline R\op \cup \Delta_{E_1})\delta =(-1)^{|E_1|}\Delta_\tau \delta \mpoint$$
Since $\tau$ permutes all the subsets $A\subseteq E$ and preserves their cardinality, we have $\Delta_\tau\delta=\delta\Delta_\tau$.
Consequently
$$\big((-1)^{|E_1|}\Delta_\tau \delta\big)^2=(-1)^{2|E_1|}\Delta_\tau^2 \delta^2=\delta^2=(-1)^{|E_1|}\Delta_\tau \delta \mvirg$$
so we obtain an idempotent.\par

Right multiplication by this idempotent defines an idempotent endomorphism of the correspondence functor~$k\CC(-,E)R$
(notice that both $\Delta_\tau$ and $\delta$ commute with~$R$).
The image of this endomorphism is the subfunctor generated by the element $(-1)^{|E_1|}\Delta_\tau \delta$, that is, 
the subfunctor generated by~$\delta$ because $\Delta_\tau$ is invertible.
But we know that ${<}\delta{>}$ is isomorphic to the fundamental functor~$\S_{E,R}$.
Therefore $\S_{E,R}$ is isomorphic to a direct summand of~$k\CC(-,E)R$, hence a direct summand of~$k\CC(-,E)$ because $R^2=R$ is idempotent.
Since $k\CC(-,E)$ is a projective functor by Yoneda's lemma, we conclude that $\S_{E,R}$ is projective.\par

Our assumption (b) also says that the $k\Aut(E,R)$-module~$V$ is projective.
By Proposition~\ref{precursor}, $S_{E,R,V}$ is isomorphic to~$\S_{E,R}\otimes_{k\Aut(E,R)} V$,
which is in turn a direct summand of~$\S_{E,R}\otimes_{k\Aut(E,R)}k\Aut(E,R) \cong \S_{E,R}$.
Therefore $S_{E,R,V}$ is projective, proving~(a).
\endpf

Another proof of the implication (b) $\Rightarrow$ (a) will be given later in Corollary~\ref{simple-SQ}.

%%%%%%%% Section

\section{Projectivity of fundamental functors} \label{Section-fundamental}

\bigskip
\noindent
Given a poset $(E,R)$, we know from Proposition~\ref{precursor} that every simple functor $S_{E,R,V}$ has a precursor $\S_{E,R}$,
called the fundamental functor associated with the poset $(E,R)$.
This is actually defined over any commutative base ring~$k$.
The main result of this section is analogous to Theorem~\ref{characterization-proj}.

\result{Theorem} \label{proj-fundamental}
Let $(E,R)$ be a finite poset. Then $\S_{E,R}$ is a projective functor if and only if $(E,R)$ is a pole poset.
\fresult

\pf
Assume first that $\S_{E,R}$ is a projective functor.
We allow the base ring $k$ to vary and we write a superscript~$(k)$ to emphasize that a functor belongs to the category $\CF_k$ of correspondence functors defined over the base ring~$k$.
Let $m$ be a maximal ideal of~$k$ and let $C=k/m$ be the corresponding field.
The scalar extension functor
$$\CF_k \longrightarrow \CF_C \;, \qquad F\mapsto C\otimes_k F$$
is left adjoint of the scalar `restriction' functor, which is obviously exact.
Therefore, scalar extension sends projective objects to projective objects.
In particular, we see that $C\otimes_k \S_{E,R}^{(k)}$ is projective.\par

By Theorem~5.6 in~\cite{BT4}, the evaluation $\S_{E,R}^{(k)}(X)$ at a finite set~$X$ has an explicit $k$-basis $\CB_X$.
This basis is defined independently of~$k$, so that it remains a $k'$-basis for any ring extension $k\to k'$.
Therefore, the natural surjection
$$C\otimes_k \S_{E,R}^{(k)}(X) \longrightarrow \S_{E,R}^{(C)}(X)$$
is an isomorphism. Since this holds for any~$X$, we have $C\otimes_k \S_{E,R}^{(k)} \cong \S_{E,R}^{(C)}$
and it follows that $\S_{E,R}^{(C)}$ is projective.
Now the functor $M:=S_{E,R}^{(C)}$ satisfies the assumptions of Lemma~\ref{proj-pole} and this lemma then asserts that $(E,R)$ is a pole poset,
as was to be shown.\par

For the converse, we use the proof of (b) $\Rightarrow$ (a) in Theorem~\ref{characterization-proj}.
This proof (except the last paragraph) tells us precisely that, whenever $(E,R)$ is a pole poset, the fundamental functor $\S_{E,R}$ is projective.
\endpf

Another proof of the projectivity of~$\S_{E,R}$ whenever $(E,R)$ is a pole poset will be given later (see Remark~\ref{proj-fundamental-2}).

%%%%%%%% Section

\section{Morphisms and idempotents corresponding to pole lattices} \label{Section-idempotents}

\bigskip
\noindent
In this section, we continue our analysis of the category~$k\CL$ of finite lattices, where $k$ is a commutative ring.
We construct morphisms involving a finite lattice~$T$ and a pole lattice~$P$.
Assuming that $\Sur_\CL(T,P)$ is nonempty, we fix a surjective join-morphism $\pi:T\to P$,
from which we will construct an idempotent endomorphism of~$T$ to the effect that $P$ can be viewed as a sort of `direct summand' of~$T$.
By means of the fully faithful functor $T\mapsto F_T$, we deduce that a certain quotient of~$F_P$ is a direct summand of~$F_T$.
We will see later in Section~\ref{Section-functors-pole} how to deduce information about projective direct summands of the correspondence functor~$F_T$.\par

Our results generalize those obtained in~\cite{BT3} in the special case of totally ordered lattices.
We follow the same line of development, but with many necessary additions and technical adaptations.\par

Recall that $P_1$ denotes the subset of elements $p\in P$ such that $p$ is comparable to every element of~$P$,
while $P_2$ denotes the subset consisting of all twins.
Let $E=\Irr(P)$ be the set of irreducible elements of~$P$, described in Lemma~\ref{irred-pole}.
We write $E_1=E\cap P_1$ and $E_2 =E\cap P_2$ (so that in fact $E_2=P_2$).\par

\result{Notation} \label{notation-B}
We define a notation associated with the surjective join-morphism $\pi:T\to P$.
\begin{enumerate}
\item For every $p\in P$, let $b_p^\pi=\pi\op(p)=\sup\big(\pi^{-1}(p)\big)$.
Whenever $\pi$ is fixed, we write simply $b_p=b_p^\pi$.
\item $B=\Im(\pi\op) = \{b_p\mid p\in P\}$.
Notice that $B$ is a subposet of~$T\op$ which is join-closed, hence a subposet of~$T$ which is meet-closed and isomorphic to~$P$.
\item For every $e\in E_1$, let $b_e^-=b_{r(e)}$ and $b_e^+=b_e$, where $r(e)=\sup[\widehat 0,e[_P$.
\item For every $e\in E_2$ and if $\breve{e}$ is the twin of~$e$, let $b_e^-=b_e$ and $b_e^+=b_{s(e)}$, where $s(e)=\inf]e,\widehat 1]_P=e\vee \breve{e}$.
\end{enumerate}
\fresult

\begin{rem}{Remark} \label{non-uniform}
The definition in (c) and~(d) is not uniform since we have $b_e=b_e^+$ in one case and $b_e=b_e^-$ in the other.
This strange behavior will be explained in Remark~\ref{uniform}, where a uniform explanation will be given.
\end{rem}

For every $e\in E$, choose $a_e\in [b_e^-,b_e^+]_T$
(where the subscript $T$ emphasizes that the interval is taken within the lattice~$T$).
This defines a family $A=(a_e)_{e\in E}$ with the following property.

\result{Lemma} \label{order-preserving}
Let $A=(a_e)_{e\in E}$ be a family of elements of~$T$ such that $a_e\in [b_e^-,b_e^+]_T$ for every $e\in E$.
Then, whenever $e,f\in E$,
$$e<_Pf \; \Longrightarrow \; a_e\leq_T a_f \mpoint$$
\fresult

\pf If $e\in E_2$ and $f\in E_1$ with $e<_P f$, then $e<_P s(e)\leq_P r(f)<_P f$ and therefore
$$a_e\leq_T b_e^+ = b_{s(e)} \leq_T b_{r(f)} =b_f^-\leq_T a_f \mpoint$$
The other cases are easier and are left to the reader.
\endpf

By Lemma~\ref{distributive-extend}, since $P$ is a distributive lattice, the order-preserving map $E\to T$, $e\mapsto a_e$, extends to a join-morphism
$$j_A^\pi:P\longrightarrow T \,,\qquad p\mapsto a_p \mpoint$$
Explicitly, we have $a_{\widehat 0}=\widehat 0$ and $a_{e\vee \breve{e}}=a_e\vee a_{\breve{e}}$ whenever $e\in E_2$ with twin~$\breve{e}$
(these are the only non-irreducible elements of~$P$ by Lemma~\ref{irred-pole}).
Note that $j_A^\pi$ is not necessarily a section of~$\pi$ (see the beginning of the proof of Proposition~\ref{rho_Y}).

Define the family $B^-=(b_e^-)_{e\in E}$ and write
$$\mu(B^-,A)=\prod_{e\in E} \mu(b_e^-,a_e)$$
where $\mu(-,-)$ denotes the M\"obius function of the lattice~$T$.
Allowing the family $A$ to vary (i.e. $a_e$ varies in $[b_e^-,b_e^+]_T$ for each~$e\in E$), define
\begin{equation}\label{def-j}
j^\pi=(-1)^{|E_1|}\sum_{A}\mu(B^-,A) \, j_A^\pi \mpoint
\end{equation}
This is a $k$-linear combination of join-morphisms, hence an element of~$k\CL(P,T)$.
The morphisms $j^\pi$ have remarkable properties, in particular when $j^\pi$ is composed with the surjection~$\pi$.
We are going to explore those properties in a series of propositions. We first start with a lemma.

\result{Lemma} \label{coincide}
Let $A=(a_e)_{e\in E}$ and $\widetilde A=(\widetilde a_e)_{e\in E}$ be two families as above and fix some $g\in E$.
Suppose that $\widetilde a_e =a_e$ for all $e\in E-\{g\}$. Then the following are equivalent~:
\begin{enumerate}
\item $j_{\widetilde A}^\pi(p)=j_A^\pi(p)$ for all $p\in P-\{g\}$.
\item If $g\in E_2$, then $\widetilde a_g\vee a_{\breve g}=a_g\vee a_{\breve g}$ where $\breve g$ is the twin of~$g$.
\end{enumerate}
\fresult

\pf 
Suppose that (b) holds. If $p=e\in E-\{g\}$, then $\widetilde a_e =a_e$ by assumption, that is, $j_{\widetilde A}^\pi(e)=j_A^\pi(e)$.
If $p=\widehat 0$, then $j_{\widetilde A}^\pi(\widehat 0)=\widehat 0=j_A^\pi(\widehat 0)$.
If now $p\in P-E$ and $p\neq\widehat 0$, then $p=u\vee \breve u$ for some $u\in E_2$, by the definition of a pole lattice.
If $g\neq u, \breve u$, then
$$j_{\widetilde A}^\pi(p)=j_{\widetilde A}^\pi(u\vee \breve u)=\widetilde a_u\vee \widetilde a_{\breve u}=a_u\vee a_{\breve u}
=j_A^\pi(u\vee \breve u)=j_A^\pi(p) \mpoint$$
If $g=u$, then the assumption~(b) implies that
$$j_{\widetilde A}^\pi(p)=j_{\widetilde A}^\pi(g\vee \breve g)=\widetilde a_g\vee \widetilde a_{\breve g}
=\widetilde a_g\vee a_{\breve g} = a_g\vee a_{\breve g} =j_A^\pi(g\vee \breve g)=j_A^\pi(p) \mvirg$$
proving (a).\mpn

Assume conversely that (a) holds. If $g\in E_1$, then condition~(b) is empty and there is nothing to prove.
So suppose that $g\in E_2$. Then
$$\widetilde a_g\vee a_{\breve g} = \widetilde a_g\vee \widetilde a_{\breve g} = j_{\widetilde A}^\pi(g\vee \breve g)
= j_A^\pi(g\vee \breve g)= a_g\vee a_{\breve g} \mvirg$$
proving (b).
\endpf

\result{Definition} \label{def-H} Associated with the subset $E=\Irr(P)$, there is a subfunctor $H_P$ of~$F_P$ defined as follows.
For any finite set~$X$, the evaluation $H_P(X)$ is the $k$-submodule of~$F_P(X)$
generated by all functions $\varphi:X\to P$ such that $E\not\subseteq \varphi(X)$.
\fresult
This subfunctor is important in the theory of correspondence functors (see Section~5 of~\cite{BT3} for details).
 
\result{Proposition} \label{vanish-on-H}
Let $\pi\in\Sur_\CL(T,P)$ and let $j^\pi:P\to T$ be the morphism defined in~(\ref{def-j}).
\begin{enumerate}
\item For any finite set $X$ and any function $\varphi:X\to P$ such that $E\not\subseteq \varphi(X)$,
we have $j^\pi\varphi=0$.
\item $j^\pi$ induces a morphism $F_{j^\pi}:F_P\to F_T$ vanishing on $H_P$, hence induces in turn a morphism
$$\overline F_{j^\pi}: F_P/H_P \longrightarrow F_T \mpoint$$
\end{enumerate}
\fresult

\pf Since (b) immediately follows from~(a), it suffices to prove~(a).
We have
$$j^\pi\varphi=(-1)^{|E_1|}\sum_A\mu(B^-,A) j_A^\pi\varphi
=\sum_{\psi:X\to T}(-1)^{|E_1|} \Big(\sum_{\substack{A \\ j_A^\pi\varphi=\psi}}\mu(B^-,A)\Big) \psi \mpoint$$
For a fixed~$\psi$, we have to prove that the inner sum over~$A$ is zero.
If this inner sum is empty, then the sum is zero and we are done.
Otherwise, we can choose $A$ such that $j_A^\pi\varphi=\psi$.
Let $g\in E$ be such that $g\notin\varphi(X)$.
Then we can modify the family~$A$ into~$\widetilde A$, by changing only the image $j_A^\pi(g)=a_g\in [b_g^-,b_g^+]_T$
into $j_{\widetilde A}^\pi(g)=\widetilde a_g\in [b_g^-,b_g^+]_T$,
with the extra condition that $\widetilde a_g\vee a_{\breve g}=a_g\vee a_{\breve g}$ in case $g\in E_2$.
The point of such a modification is that it is precisely the only kind
which does not change the equality $j_A^\pi\varphi=\psi$, by Lemma~\ref{coincide}.
We set $A'=(a_e)_{e\in E-\{g\}}$ and $B'^-=(b_e^-)_{e\in E-\{g\}}$ and we let
$$j_{A'}^\pi: E-\{g\}  \longrightarrow T\,, \qquad e\mapsto a_e \mvirg$$
which we extend to a join morphism $j_{A'}^\pi:P-\{g\}\to T$.
We obtain
$$\sum_{\substack{A \\ j_A^\pi\varphi=\psi}}\mu(B^-,A)
= \sum_{\substack{A' \\ j_{A'}^\pi\varphi=\psi}}\mu(B'^-,A') \sum_{\widetilde a_g} \mu(b_g^-,\widetilde a_g) \mvirg$$
where the inner sum runs over all $\widetilde a_g\in [b_g^-,b_g^+]_T$,
with the extra condition that $\widetilde a_g\vee a_{\breve g}=a_g\vee a_{\breve g}$ in case $g\in E_2$.\par

If $g\notin E_2$, then the sum runs over all $\widetilde a_g\in [b_g^-,b_g^+]_T$
and this is zero by the definition of the M\"obius function (because $b_g^-=b_{r(g)}<_T b_g= b_g^+$).
If $g\in E_2$, then the extra condition is equivalent to $\widetilde a_g\vee (b_g\vee a_{\breve g})=a_g\vee a_{\breve g}$
(because $b_g=b_g^-\leq \widetilde a_g$),
so $\widetilde a_g$ runs over the interval $[b_g,a_g\vee a_{\breve g}]_T$ with the condition that
its join with the fixed element $b_g\vee a_{\breve g}$ is equal to the top element $a_g\vee a_{\breve g}$.
By a well-known property of the M\" obius function (Corollary~3.9.3 in \cite{St}), the corresponding sum
$$\sum_{\substack{\widetilde a_g\in [b_g,a_g\vee a_{\breve g}]_T \\
\widetilde a_g\vee (b_g\vee a_{\breve g})=a_g\vee a_{\breve g}}}
\mu(b_g^-,\widetilde a_g)$$
is zero, provided the fixed element $b_g\vee a_{\breve g}$ is not equal to the bottom element~$b_g$.
But this is indeed the case since $b_g\vee a_{\breve g}\geq_T b_g\vee b_{\breve g} >_T b_g$,
the latter inequality coming from the fact that $\pi(b_g\vee b_{\breve g})=g\vee\breve g >_P g=\pi(b_g)$.
It follows that the coefficient of every $\psi$ is zero, hence $j^\pi\varphi=0$.\endpf

Now we want to compute the composite $\pi j^\pi$. For any subset $Y$ of~$E$, we define
\begin{equation} \label{def-rho-Y}
\rho_Y: E \longrightarrow P \,, \qquad \rho_Y(e)=
\begin{cases}
{e} & {\text { if }\, e\in Y \mvirg} \\
{r(e)} & {\text { if }\, e\in E_1, \, e\notin Y \mvirg} \\
{s(e)} & {\text { if }\, e\in E_2, \, e\notin Y \mpoint} \end{cases}
\end{equation}
It is easy to see that $\rho_Y$ is order-preserving (because, if $e\in E_2$, $f\in E_1$, and $e<_P f$, then $e<_P s(e)\leq_P r(f) <_P f$,
while the other cases are easier). Therefore, by Lemma~\ref{distributive-extend}, $\rho_Y$ extends to a join-morphism $\rho_Y: P \longrightarrow P$
because the pole lattice~$P$ is distributive.
Note that $\rho_Y(p)=p$ for any $p\notin E$. This is clear if $p=\widehat 0$.
Otherwise $p=e\vee \breve e$ for some $e\in E_2$ by Lemma~\ref{irred-pole} and
$$p=e\vee \breve e\leq_P \rho_Y(e)\vee \rho_Y(\breve e)\leq_P s(e)\vee s(\breve e)=p\vee p=p \mvirg$$
forcing equality and $\rho_Y(p)=\rho_Y(e)\vee \rho_Y(\breve e)=p$.

\result{Proposition} \label{rho_Y}
Let $\pi\in\Sur_\CL(T,P)$ and let $j^\pi:P\to T$ be the morphism defined in~(\ref{def-j}).
\begin{enumerate}
\item $\pi j^\pi= \displaystyle \sum_{\emptyset \subseteq Y \subseteq E} (-1)^{|E-Y|} \rho_Y$, where $\rho_Y$ is defined by~(\ref{def-rho-Y}).
\item If $Y\neq E$, then $E\not\subseteq \rho_Y(P)$.
\end{enumerate}
\fresult

\pf
For simplicity, we write $<$ instead of ~$<_P$ and $\leq$ instead of~$\leq_P$.\par

(a) If $e\in E_1$ and $b_{r(e)}<x$ in~$T$, then $r(e)< \pi(x)$ because $b_{r(e)}=\sup\{t\in T\mid \pi(t)=r(e)\}$.
Thus if $x\in \, ]b_{r(e)},b_e]_T$, we get $r(e)< \pi(x)\leq e$, hence $\pi(x)=e$.
Similarly, if $e\in E_2$ and $x\in \, ]b_e,b_{s(e)}]_T$, then $\pi(x)=s(e)$.
It follows that
$$\pi j_A^\pi(e)=\begin{cases} {e}& {\text { if }\, e\in E_1 \text{ and } j_A^\pi(e)\in \, ]b_e^-,b_e^+]_T = ]b_{r(e)},b_e]_T \mvirg} \\
{e}& {\text { if }\, e\in E_2 \text{ and } j_A^\pi(e)=b_e^-=b_e \mvirg} \\
{r(e)}& {\text { if }\, e\in E_1 \text{ and } j_A^\pi(e)=b_e^-=b_{r(e)} \mvirg} \\
{s(e)}& {\text { if }\, e\in E_2 \text{ and } j_A^\pi(e)\in \, ]b_e^-,b_e^+]_T = ]b_e,b_{s(e)}]_T \mpoint} \end{cases}$$
We see that $\pi j_A^\pi=\rho_Y$ for a suitable subset $Y\subseteq E$ and therefore
$$\pi j^\pi=\sum_{\emptyset \subseteq Y \subseteq E} (-1)^{|E_1|} 
\Big( \sum_{\substack{A \; \\ \pi j_A^\pi=\rho_Y}}\mu(B^-,A) \Big) \rho_Y \mpoint$$
For a fixed subset $Y$, in order to realize the condition $\pi j_A^\pi=\rho_Y$, we have the following possibilities~:
\begin{enumerate}
\item[$\bullet$] If $e\in Y\cap E_1$, then $j_A^\pi(e)$ can run freely in~$]b_e^-,b_e^+]_T$.
\item[$\bullet$] If $e\in Y\cap E_2$, then $j_A^\pi(e)$ must be equal to~$b_e^-=b_e$.
\item[$\bullet$] If $e\in (E-Y)\cap E_1$, then $j_A^\pi(e)$ must be equal to~$b_e^-=b_{r(e)}$.
\item[$\bullet$] If $e\in (E-Y)\cap E_2$, then $j_A^\pi(e)$ can run freely in~$]b_e^-,b_e^+]_T$.
\end{enumerate}
It follows that the coefficient $\displaystyle(-1)^{|E_1|} \sum_{\substack{A \; \\ \pi j_A^\pi=\rho_Y}} \mu(B^-,A)$ is equal to
$$
\begin{array}{l}
\displaystyle(-1)^{|E_1|} \prod_{e\in Y\cap E_1} \Big( \sum_{a_e\in \,]b_e^-,b_e^+]_T} \mu(b_e^-,a_e) \Big)
\cdot \prod_{e\in (E-Y)\cap E_2}  \Big( \sum_{a_e\in \,]b_e^-,b_e^+]_T} \mu(b_e^-,a_e) \Big) \\
= (-1)^{|E_1|} \cdot (-1)^{|Y\cap E_1|} \cdot (-1)^{|(E-Y)\cap E_2|} \\
= (-1)^{|(E-Y)\cap E_1|}  \cdot (-1)^{|(E-Y)\cap E_2|} \\
= (-1)^{|E-Y|}  \mvirg
\end{array}$$
using the fact that
$$0=\sum_{a_e\in \,[b_e^-,b_e^+]_T} \mu(b_e^-,a_e) =1+ \sum_{a_e\in \,]b_e^-,b_e^+]_T} \mu(b_e^-,a_e) \mpoint$$
This shows that
$$\pi j^\pi=\sum_{\emptyset \subseteq Y \subseteq E} (-1)^{|E-Y|} \rho_Y \mvirg$$
as required.\mpn

(b) Suppose that $Y$ is a proper subset of~$E$ and let $g\in E$ be maximal such that $g\notin Y$.
We want to prove that $g\notin \rho_Y(P)$. We let $p\in P$ and we prove that $\rho_Y(p)\neq g$.\par

If $p>g$, then $p\in Y$ and $\rho_Y(p)=p\neq g$, while if $p=g$, then $\rho_Y(g)\neq g$.\par

Assume now that $p\not\geq g$ and $g\in E_1$. Then $p<g$ and $\rho_Y(p)\leq \rho_Y(g)=r(g)<g$.\par

Assume now that $p\not\geq g$ and $g\in E_2$. Then either $p<g$ or $p=\breve g$, the twin of~$g$.

If $p< g$ and $p\in E_1$, then $\rho_Y(p)\leq p<g$.
If $p< g$ and $p\in E_2$, then $\rho_Y(p)\leq s(p)\leq g$.
But $s(p)$ is reducible since $s(p)=p\vee \breve p$, while $g$ is irreducible. Therefore $s(p)\neq g$, hence $\rho_Y(p)<g$.\par

If $p=\breve g$, then $p\in E_2$ and $\rho_Y(p)$ is either $p$ or $s(p)$. But neither $p$ nor $s(p)$ is equal to~$g$.\par

We have proved that $\rho_Y(p)\neq g$ in all cases, as was to be shown.
\endpf

\result{Remark} \label{uniform}
{\rm 
In the special case when $T=P$ and $\pi=\Id$, we find that $j^{\Id}$ is a linear combination of the maps~$\rho_Y$.
It turns out that $j^{\Id}$ is actually an avatar of the element $\gamma_{P\op}\in F_P(E^0)$ which is defined in~(\ref{def-gamma}), where $E^0=\Irr(P\op)$.
We know that the element $\gamma_T$ plays an important role throughout the theory of correspondence functors
(see Section~9 of~\cite{BT3} and Section~\ref{Section-characterization} of the present paper).
The advantage of $\gamma_{P\op}$ is that it has a uniform definition, contrary to~$j^P$ (as observed in Remark~\ref{non-uniform}).\par

To make this explicit, let $E=\Irr(P)$, viewed as a subposet of~$P$ and $E^0=\Irr(P\op)$, viewed also as a subposet of~$P$
(so that it is actually $(E^0)\op$ which is the subposet of irreducible elements of~$P\op$).
Since $P$ is a distributive lattice, it is isomorphic to the lattice $\Idown(E)$ of all subsets of~$E$ closed under taking smaller elements.
The passage to complements induces an isomorphism $\Idown(E)\cong \Iup(E)\op$, where $\Iup(E)$ is the lattice of all subsets of~$E$ closed under taking greater elements. On restriction to~$E$, this induces an order-preserving isomorphism $\alpha : E \to E^0$, which turns out to map $e\in E_1$ to~$r(e)\in E^0_1$ (in the totally ordered part) and $e\in E_2$ to its twin~$\breve e\in E^0_2$ (in the twin part).\par

Now $\gamma_{P\op}$ is a linear combination of maps $E^0\to P$ and we precompose it with $\alpha\tau$,
where $\tau:E\to E$ exchanges all the twins and fixes all the other elements.
We obtain a linear combination of maps $E\to P$ and, after an explicit computation, it turns out that
$$\gamma_{P\op} \, \alpha\,\tau = \pm\,  j^{\Id} \mvirg$$
the sign being actually $(-1)^{|E_1|}$.
(This computation appears explicitly in the proof of Theorem~\ref{S_Q-fundamental},
using a bijection $\omega:E^0\to E$ which is actually the inverse of~$\alpha\tau$.)
The definition of~$\rho_Y$ in~(\ref{def-rho-Y}) was {\em not} uniform and, accordingly, $j^{\Id}$ has a rather strange behavior.
However, by means of the isomorphism $\alpha\tau$, the translation of all this in terms of~$\gamma_{P\op}$ becomes uniform.\par

Unfortunately, we need to work with $j^{\Id}$ rather than $\gamma_{P\op}$.
The reason is that $\gamma_{P\op}\in F_P(E^0)$ is a linear combination of maps $E^0\to P$,
whereas, after composing with~$\alpha$, we obtain order-preserving maps $E\to P$ which are extendible to endomorphisms $P\to P$
(because $P$ is a distributive lattice, see Lemma~\ref{distributive-extend}).
The key fact is that endomorphisms are better because they can be composed, in particular it makes sense to consider idempotent endomorphisms.
}
\fresult

\bigskip
We can now prove a main result concerning the composite $j^\pi\pi$
and obtain consequences for the correspondence functor $F_P$ associated with the pole lattice~$P$.

\result{Proposition} \label{ker-j}
Let $\pi\in\Sur_\CL(T,P)$ and let $j^\pi:P\to T$ be the morphism defined in~(\ref{def-j}).
Let $q:F_P \to F_P/H_P$ be the canonical surjection, where $H_P$ is defined by~(\ref{def-H}).
\begin{enumerate}
\item $j^\pi\pi$ is an idempotent endomorphism of~$T$.
\item The composite of $\overline F_{j^\pi}: F_P/H_P \to F_T$ and $q\, F_{\pi}: F_T\to F_P/H_P$ is the identity morphism of~$F_P/H_P$.
\item $\overline F_{j^\pi}:F_P/H_P \to F_T$ is injective and embeds $F_P/H_P$ as a direct summand of~$F_T$.
\item $F_{j^\pi}F_{\pi}$ is an idempotent endomorphism of~$F_T$ whose image is isomorphic to~$F_P/H_P$.
\end{enumerate}
\fresult

\pf
(a) This follows from~(d), which is proved below, because the functor $F_?:k\CL \to \CF_k$ is fully faithful by Theorem~\ref{fully-faithful}.
Alternatively, it is not difficult to compute directly
$$j^\pi\pi j^\pi=j^\pi \sum_{\emptyset \subseteq Y \subseteq E} (-1)^{|E-Y|} \rho_Y
=j^\pi \Id_P + \sum_{ Y \neq E} (-1)^{|E-Y|} j^\pi\rho_Y=j^\pi \mvirg$$
because $E\not\subseteq\rho_Y(P)$ if $Y\neq E$ by Proposition~\ref{rho_Y}, hence $j^\pi\rho_Y=0$ by Proposition~\ref{vanish-on-H}.
Then the equality $j^\pi\pi j^\pi=j^\pi$ implies that $j^\pi\pi$ is an idempotent.\mpn

(b) By Proposition~\ref{rho_Y}, for any finite set~$X$ and any function $\varphi : X\to P$, 
$$F_{\pi} F_{j^\pi} (\varphi) = \pi j^\pi \varphi = \sum_{\emptyset \subseteq Y \subseteq E} (-1)^{|E-Y|} \rho_Y\varphi
=\varphi + \sum_{Y\neq E} (-1)^{|E-Y|} \rho_Y\varphi \mpoint$$
But $E\not\subseteq\rho_Y(P)$ if $Y\neq E$ by Proposition~\ref{rho_Y}, hence $E\not\subseteq\rho_Y\varphi(X)$.
In other words, $\rho_Y\varphi\in H_P(X)$, so that
$$F_{\pi} F_{j^\pi} (\varphi) = \varphi \pmod{H_P(X)} \mpoint$$
Composing with the canonical map $q:F_P(X) \to F_P(X)/H_P(X)$ and writing $q(\varphi)=\overline\varphi$, we obtain
$$q F_{\pi}\overline F_{j^\pi} (\overline\varphi) = q F_{\pi} F_{j^\pi}(\varphi) = q(\varphi)= \overline\varphi \mvirg$$
as was to be shown.\mpn

(c) This follows immediately from~(b).\mpn

(d) This follows immediately from~(b) and the obvious equality $\overline F_{j^\pi}q\,F_{\pi} = F_{j^\pi}F_{\pi}$.\mpn
\endpf

One of our aims is to show that the idempotents $j^\pi\pi$ are orthogonal.
In order to understand the product of two idempotents $j^\theta\theta$ and $j^\pi\pi$ we need to have more information about $\theta j^\pi$.
This is the purpose of our next three propositions, but we first need a lemma.

\result{Lemma} \label{not-singleton}
Let $Q$ be a pole lattice, let $\theta\in\Sur_\CL(T,Q)$, and let $[t_1,t_2]_T$ be an interval in~$T$.
For every $q\in Q$, define
$$U^q=\theta^{-1}(q) \cap [t_1,t_2]_T = \{ a\in [t_1,t_2]_T \mid \theta(a)=q)\} \subseteq [t_1,t_2]_T \mpoint$$
Let $q_1=\theta(t_1)$ and assume that $U^{q_1}$ is not reduced to the singleton~$t_1$.
Then for each $q\in Q$, we have $\displaystyle\sum_{a\in U^q} \mu(t_1,a)=0$.
\fresult

\pf
The result is obvious if $U^q=\emptyset$, so we assume that $U^q$ is nonempty.
Since all elements of~$U^q$ have the same image under~$\theta$, so has their join and therefore $U^q$ has a supremum
$$u^q=\sup (U^q) \in U^q \mpoint$$
Now we have $[t_1,t_2]_T=\sqcup_{q\in Q} U^q$
and, by assumption, $U^{q_1} =[t_1,u^{q_1}]_T$ is a nontrivial interval, so that
$$\sum_{a\in U^{q_1}} \mu(t_1,a)=0 \mpoint$$
This is the starting point of an induction argument.
We fix $q\in Q$ and we assume by induction that $\displaystyle\sum_{a\in U^r} \mu(t_1,a)=0$ for every $r\in Q$ such that $q_1\leq r < q$.
Then we obtain
$$\begin{array}{rcl}
0&=&\displaystyle\sum_{a\in[t_1,u^q]_T}\mu(t_1,a)= \sum_{q_1\leq r \leq q} \, \sum_{a\in U^r} \mu(t_1,a) \\
&=& \displaystyle \sum_{a\in U^q} \mu(t_1,a) + \sum_{q_1\leq r < q} \, \sum_{a\in U^r} \mu(t_1,a) \\
&=&\displaystyle\sum_{a\in U^q} \mu(t_1,a) \mvirg
\end{array}$$
using the induction assumption. This completes the proof.
\endpf

\result{Proposition} \label{compose-nonzero}
Let $\pi\in\Sur_\CL(T,P)$ and $\theta\in\Sur_\CL(T,Q)$, where $P$ and $Q$ are pole lattices,
and let $j^\pi:P\to T$ be the morphism defined in~(\ref{def-j}).
Suppose that $\theta j^\pi\neq 0$. Then the restriction of $\theta$ to the subset $B=\Im(\pi\op)$ is injective.
In particular, $|P|\leq |Q|$.
\fresult

\pf Let $E=\Irr(P)$. By the definition of~$j^\pi$, we have
$$\theta j^\pi=(-1)^{|E_1|}\sum_A\mu(B^-,A) \, \theta j_A^\pi
= \sum_{\psi:P\to Q} \Big((-1)^{|E_1|} \sum_{\substack{A\; \\ \theta j_A^\pi=\psi}}\mu(B^-,A)\Big) \psi \mpoint$$
Now fix some morphism $\psi:P\to Q$ and, for every $e\in E$ and every $q\in Q$, define
$$U_e^q=\theta^{-1}(q) \cap [b_e^-,b_e^+]_T = \{ a\in [b_e^-,b_e^+]_T \mid \theta(a)=q)\} \subseteq [b_e^-,b_e^+]_T \mpoint$$
Here, we write $B=\{b_p\mid p\in P\}$, as before.
Then, since a join-morphism from $P$ is entirely determined on $E=\Irr(P)$, we have
$$\theta j_A^\pi=\psi \iff j_A^\pi(e)\in U_e^{\psi(e)} \;\forall e\in E \iff a_e\in U_e^{\psi(e)} \;\forall e\in E \mpoint$$
In particular, if $\psi$ appears in the expression of~$\theta j^\pi$, then $U_e^{\psi(e)}\neq\emptyset$ for every $e\in E$.
It follows now that the coefficient of~$\psi$ is, up to sign, equal to
$$\sum_{\substack{A\; \\ \theta j_A^\pi=\psi}}\mu(B^-,A) = \prod_{e\in E} \Big( \sum_{a_e\in U_e^{\psi(e)}} \mu(b_e^-,a_e) \Big) \mpoint$$
Suppose that $\theta_{|B} : B\to Q$ is not injective. Then we want to prove that the coefficient of~$\psi$ is zero.
This is the case if $U_e^{\psi(e)} = \emptyset$ for some $e\in E$, because we get an empty sum, which is zero.
So we assume that $U_e^{\psi(e)} \neq \emptyset$ for every $e\in E$.
The noninjectivity of $\theta_{|B}$ implies that there exist two adjacent elements $w<y$ in~$P$ such that $\theta(b_w)=\theta(b_y)$.
There are three cases.\mpn

{\bf Case 1.} $y\in E_1$ and $w=r(y)$. Then $b_w=b_y^-$ and $b_y=b_y^+$.
Choosing $a\in U_y^{\psi(y)}$, we obtain
$$\theta(b_w)=\theta(b_y^-)\leq_Q \theta(a) \leq_Q \theta(b_y^+) = \theta(b_y) \mvirg$$
hence $\theta(b_y^-)=\theta(a)=\theta(b_y^+)$.
Since $\theta(a)=\psi(y)$, it follows that the whole interval $[b_y^-,b_y^+]_T$ is mapped to~$\psi(y)$ under~$\theta$, that is,
$[b_y^-,b_y^+]_T=U_y^{\psi(y)}$. But then
$$\sum_{a_y\in U_y^{\psi(y)}} \mu(b_y^-,a_y) =0 \mvirg$$
by the definition of the M\"obius function (because $b_y^-\neq b_y^+$). Therefore the coefficient of~$\psi$ is zero.\mpn

{\bf Case 2.} $w\in E_2$ and $y=s(w)$. Then $b_w=b_w^-$ and $b_y=b_w^+$.
Choosing $a\in U_w^{\psi(w)}$, we obtain
$$\theta(b_w)=\theta(b_w^-)\leq_Q \theta(a) \leq_Q \theta(b_w^+) = \theta(b_y) \mvirg$$
hence $\theta(b_w^-)=\theta(a)=\theta(b_w^+)$.
Since $\theta(a)=\psi(w)$, it follows that the whole interval $[b_w^-,b_w^+]_T$ is mapped to~$\psi(w)$ under~$\theta$, that is,
$[b_w^-,b_w^+]_T=U_w^{\psi(w)}$. But then
$$\sum_{a_w\in U_w^{\psi(w)}} \mu(b_w^-,a_w) =0 \mvirg$$
and the coefficient of~$\psi$ is zero.\mpn

{\bf Case 3.} $y\in E_2$ and $w=r(y)$. Let $z=\breve y$ be the twin of~$y$, so that $b_z^-=b_z$. Since $w<z$, we have $b_w< b_z$ and
$$\theta(b_y\vee b_z)=\theta(b_y)\vee \theta(b_z)=\theta(b_w)\vee \theta(b_z)=\theta(b_w\vee b_z)=\theta(b_z) \mpoint$$
Letting $q_1=\theta(b_z)$, we see that $U_z^{q_1}$ contains both the minimal element $b_z=b_z^-$ of the interval $[b_z^-,b_z^+]_T$
and another element $b_y\vee b_z$, because $b_z<b_y\vee b_z\leq b_{y\vee z}=b_z^+$.
Thus the assumption of Lemma~\ref{not-singleton} is satisfied and it follows that
$$\sum_{a_z\in U_z^{\psi(z)}} \mu(b_z^-,a_z) =0 \mpoint$$
Again the coefficient of~$\psi$ is zero and we are done.\par

This completes the proof of the injectivity of $\theta_{|B} : B\to Q$.
Since $\pi\op$ is injective (by Lemma~\ref{op-morphisms}), its image~$B$ has cardinality~$|P|$ and therefore $|P|\leq |Q|$.
\endpf

\result{Proposition} \label{unique-iso}
Let $\pi\in\Sur_\CL(T,P)$ and $\theta,\chi\in\Sur_\CL(T,Q)$, where $P$ and $Q$ are pole lattices.
Suppose that $j^\chi\theta j^\pi\neq 0$.
\begin{enumerate}
\item There exists a unique isomorphism $\tau: P \to Q$ of lattices such that
$$\theta j^\pi=\tau \pmod{\Ker(j^\chi)} \mpoint$$
\item Moreover, $\theta(b_p)=\tau(p)$, for all $p\in P$ (where $b_p=b^\pi_p=\pi\op(p)$, as before).
\end{enumerate}
\fresult

\pf
We assume that $j^\chi\theta j^\pi\neq 0$, and in particular $\theta j^\pi\neq 0$. Let $E=\Irr(P)$ and write first
$$\theta j^\pi=(-1)^{|E_1|}\sum_A\mu(B^-,A) \, \theta j_A^\pi
= \sum_{\psi:P\to Q} \Big((-1)^{|E_1|} \sum_{\substack{A\; \\ \theta j_A^\pi=\psi}}\mu(B^-,A)\Big) \psi \mpoint$$
Let $\psi:P\to Q$ be a map appearing with a nonzero coefficient in the expression of~$\theta j^\pi$ and
let $A$ be such that $\theta j_A^\pi=\psi$. 
Since $j^\chi\theta j^\pi\neq 0$, we can also assume that $\psi$ is such that $j^\chi\psi\neq0$.
Proposition~\ref{vanish-on-H} implies that the function $\psi:P\to Q$ must satisfy $\Irr(Q)\subseteq \psi(P)$.
Since $\psi=\theta j_A^\pi$ is a join-morphism and $\Irr(Q)$ generates~$Q$, the map $\psi:P\to Q$ must be surjective.
By Proposition~\ref{compose-nonzero}, $\theta j^\pi\neq 0$ implies that $|P|\leq |Q|$.
Therefore $|P|= |Q|$. It follows that $\psi$ is a bijective join-morphism, hence an isomorphism of lattices.\par

Proposition~\ref{compose-nonzero} also asserts that the map $\theta_{|B} : B\to Q$ is injective.
Since $|B|=|P|=|Q|$, it is a bijection and therefore there is a unique isomorphism $\tau:P\to Q$ such that
$$\theta(b_p)=\tau(p) \,, \qquad \forall \, p\in P \mpoint$$
For any $e\in E$, we have $\psi(e)=\theta j_A^\pi(e)=\theta(a_e)$ for some $a_e \in [b_e^-,b_e^+]_T$.
If $e\in E_1$, then $b_e^+=b_e$, hence
$$\psi(e)=\theta(a_e) \leq_Q \theta(b_e)=\tau(e) \mpoint$$
Therefore $e\leq_P \psi^{-1}\tau(e)$, so that $\psi^{-1}\tau(e)=e$ because $\psi^{-1}\tau$ is an automorphism of~$P$, hence height-preserving.
Similarly, if $e\in E_2$, then $b_e^-=b_e$, hence
$$\tau(e)=\theta(b_e) \leq_Q \theta(a_e)=\psi(e) \mvirg$$
so that $\psi^{-1}\tau(e)\leq_P e$ and $\psi^{-1}\tau(e)=e$.
This shows that $\psi_{|E}=\tau_{|E}$, hence $\psi=\tau$.
Therefore, whenever $A$ is such that $j^\chi\theta j_A^\pi\neq 0$, then $\theta j_A^\pi=\tau$.
It follows that the functions $\psi$ which appear with a nonzero coefficient in the expression of~$\theta j^\pi$
are $\tau$ and maps in the kernel of~$j^\chi$.
\endpf

In the situation of Proposition~\ref{unique-iso}, we can replace $\theta$ by $\theta':=\tau^{-1}\theta$ and $j^\chi$ by $j^{\chi'}:=j^\chi\,\tau$.
The effect of this is that we are reduced to the case where $Q=P$ and $\tau=\Id_P$, that is,
$$\theta(b_p)=p \,, \qquad \forall \, p\in P \mpoint$$
For simplicity, we use this reduction in our final result, which is the key for understanding the composition of the morphisms we have introduced.

\result{Proposition} \label{compute-composite}
Let $\pi,\theta,\chi\in\Sur_\CL(T,P)$, where $P$ is a pole lattice.
Suppose that $\theta(b_p)=p$ for all $p\in P$ (where $b_p=b^\pi_p=\pi\op(p)$, as before).
\begin{enumerate}
\item If $j^\chi\theta j^\pi\neq 0$, then $\theta=\pi$.
\item We have
$$j^\chi\theta j^\pi=\begin{cases}{j^\chi} & { \text{ if } \, \theta=\pi \,,} \\
{0} & { \text{ if } \, \theta\neq \pi \,.}\end{cases}$$
\end{enumerate}
\fresult

\pf
(a) By Proposition~\ref{unique-iso}, we have $\theta j^\pi=\Id_P \pmod{\Ker(j^\chi)}$, because the automorphism $\tau$ is the identity by assumption.
Moreover, as in the proof of the previous propositions, the coefficient of~$\Id_P$ in the expression of $\theta j^\pi$ is equal to
$$(-1)^{|E_1|} \sum_{\substack{A\; \\ \theta j_A^\pi=\Id_P}}\mu(B^-,A)
=(-1)^{|E_1|} \prod_{e\in E} \Big( \sum_{a_e\in U_e} \mu(b_e^-,a_e) \Big) \mvirg$$
where we write simply
$$U_e:=U_e^e= \{ a\in [b_e^-,b_e^+]_T \mid \theta(a)=e\} \subseteq [b_e^-,b_e^+]_T \mpoint$$
Since the coefficient of $\Id_P$ in the expression of $\theta j^\pi$ is nonzero (it is~1), every sum $\sum_{a_e\in U_e} \mu(b_e^-,a_e)$ is nonzero,
and in particular $U_e\neq\emptyset$.\par

As in the proof of Lemma~\ref{not-singleton}, $U_e$ has a supremum~$u_e\in U_e$.
We also define
$$V_e:= \{ a\in [b_e^-,b_e^+]_T \mid \theta(a)\leq_P r(e)\}$$
so that $[b_e^-,u_e]_T = V_e\sqcup U_e$, because any $a\leq_T u_e$ satisfies $\theta(a)\leq_P e$, hence either $\theta(a)\leq r(e)$ or $\theta(a)=e$.
There are two cases.\mpn

{\bf Case A.} $V_e=\emptyset$. Then $U_e=[b_e^-,u_e]_T$. The nonzero sum $\sum_{a_e\in U_e} \mu(b_e^-,a_e)$ forces $b_e^-=u_e$, hence
$U_e=\{b_e^-\}$.\mpn

{\bf Case B.} $V_e\neq\emptyset$. Then again $V_e$ must have a supremum $v_e\in V_e$, so that $V_e=[b_e^-,v_e]_T$ and $v_e<_Tu_e$.
In that case, we obtain
$$0\neq \sum_{a_e\in U_e} \mu(b_e^-,a_e) = \sum_{a_e\in [b_e^-,u_e]_T} \mu(b_e^-,a_e) \; - \sum_{a_e\in V_e} \mu(b_e^-,a_e)
= -\sum_{a_e\in V_e} \mu(b_e^-,a_e) \mvirg$$
because the sum over $ [b_e^-,u_e]_T$ is zero since $b_e^- \leq_T v_e <_T u_e$.
Therefore
$$\sum_{a_e\in [b_e^-,v_e]_T} \mu(b_e^-,a_e)=\sum_{a_e\in V_e} \mu(b_e^-,a_e)\neq 0$$
and this forces $b_e^-=v_e$, hence $U_e=]b_e^-,u_e]$.\mpn

By assumption, we know that $\theta(b_e)=e$ for all $e\in E$, hence $b_e\in U_e$.
If $e\in E_2$, then $b_e=b_e^-$, hence $b_e^-\in U_e$. This forces to be in case~$A$ and therefore we obtain~:\mpn

{\bf Case A.} $U_e=\{b_e\}$ if $e\in E_2$.\mpn

If $e\in E_1$, then $b_e=b_e^+$, hence $b_e^+\in U_e$. This forces to be in case~$B$ with moreover $u_e=b_e^+$.
Since $b_e^-=b_{r(e)}$, we get~:\mpn

{\bf Case B.} $U_e=]b_{r(e)},b_e]$ if $e\in E_1$.\mpn

Let $c_p=\theta\op(p)=\sup\big(\theta^{-1}(p)\big)$ (that is, $c_p=b_p^\theta$ using Notation~\ref{notation-B}).
Since we assume that $\theta(b_p)=p$, we have $b_p\leq_T c_p$ for all $p\in P$.
We now prove that $b_p=c_p$ by descending induction in the lattice~$T$, starting from the obvious equality $b_{\widehat 1_P}=\widehat 1_T =c_{\widehat 1_P}$. For simplicity, we write $<$ and $\leq$ for the order relation in~$T$.
Suppose now that $p\in P$ and $b_q=c_q$ for every $q>p$. We have to discuss three cases.\par

Assume that $p=r(e)$ with $e\in E_1$. Then $b_p\leq c_p< c_e=b_e$, hence $c_p\in [b_{r(e)},b_e]_T=\{b_p\}\sqcup U_e$ (Case~B).
But $\theta(c_p)=p\neq e$, so $c_p\notin U_e$. Therefore $c_p=b_p$.\par

Assume that $p\in E_2$. Then
$$b_p^-=b_p\leq c_p < c_{s(p)}=b_{s(p)}=b_p^+ \mpoint$$
Therefore $c_p\in [b_p^-,b_p^+]_T = U_p\sqcup ]b_p^-,b_p^+]_T$ (Case~A). Since $\theta(c_p)=p$, we have $c_p\in U_p=\{b_p\}$, hence $c_p=b_p$.\par

Assume now that $p=e\wedge \breve e$ where $e\in E_2$ with twin $\breve e$. Then $b_e=c_e$ and $b_{\breve e}=c_{\breve e}$.
Thus we obtain
$$b_p=b_{e\wedge \breve e}=b_e\wedge b_{\breve e}=c_e\wedge c_{\breve e}=c_{e\wedge \breve e}=c_p \mvirg$$
as was to be shown. We have now covered all cases, completing the proof that $b_p=c_p$ for all $p\in P$.\par

Now we obtain $\theta\op(p)=c_p=b_p=\pi\op(p)$ for all $p\in P$, hence $\theta\op=\pi\op$.
Passing to the opposite, it follows that $\theta=\pi$, as was to be shown.
\mpn

(b) We now know by~(a) that $\theta=\pi$ whenever $j^\chi\theta j^\pi\neq 0$.
Moreover, in that case, Proposition~\ref{unique-iso} implies that $\pi j^\pi$ is the sum of~$\Id_P$ and morphisms in the kernel of~$j^\chi$,
using our assumption that $\theta(b_p)=p$, for all $p\in P$.
Applying $j^\chi$, it follows that $j^\chi\pi j^\pi=j^\chi$.
\endpf

Keeping our fixed finite lattice~$T$, we now allow the pole lattice $P$ to vary.

\result{Notation} \label{representatives}
\begin{enumerate}
\item $\Pol_T$ is a set of representatives of isomorphism classes of pole lattices~$P$ such that $\Sur_\CL(T,P)$ is nonempty (hence in particular $|P|\leq |T|$, so that $\Pol_T$ is finite).
\item For any $P\in \Pol_T$, the group $\Aut(P)$ acts on $\Sur_\CL(T,P)$ (by composition)
and we let $\overline{\Sur_\CL(T,P)}$ be a fixed chosen set of representatives of the orbits.
\item If $\chi,\theta\in \overline{\Sur_\CL(T,P)}$ and $\tau\in\Aut(P)$, we define
$$f_{\chi,\tau,\theta}=j^\chi\tau\theta : T \longrightarrow T \mpoint$$
In particular, $f_{\pi,\Id_P,\pi}=j^\pi\pi$ is the idempotent of Proposition~\ref{ker-j}.
\end{enumerate}
\fresult

\result{Remark} \label{Aut-action}
{\rm
\begin{enumerate}
\item Let $\chi'=\sigma\chi\in \Sur_\CL(T,P)$ be the image of~$\chi$ under the action of~$\sigma$, for some $\sigma\in\Aut(P)$.
Then $j^{\chi'}=j^\chi\sigma^{-1}$.
This is proved by going back to Notation~\ref{notation-B} and using the associated elements
$b^\chi_p=\chi\op(p)$, respectively
$$b^{\chi'}_p={\chi'}\op(p)=\chi\op(\sigma\op(p))=\chi\op\sigma^{-1}(p)=b^\chi_{\sigma^{-1}(p)} \mvirg$$
from which the associated morphism $j^\chi$, respectively $j^{\chi'}$, is constructed, as in~(\ref{def-j}).
It is then elementary to check that $j^{\chi'}=j^\chi\sigma^{-1}$.
\item 
Changing the choice of orbit representatives has the following effect.
Let $\chi'=\sigma\chi\in \Sur_\CL(T,P)$ and $\theta'=\rho\theta\in \Sur_\CL(T,P)$, where $\sigma, \rho\in\Aut(P)$.
It follows from (a) that we obtain
$j^{\chi'} \sigma\tau\rho^{-1}\theta'=j^\chi\tau\theta$.
\item In particular, $f_{\pi,\Id_P,\pi}=j^\pi\pi$ is independent of the choice of $\pi$ in its $\Aut(P)$-orbit.
\end{enumerate}
} \fresult

Now we come to the crucial relations among the endomorphisms $f_{\chi,\tau,\theta}$.

\result{Theorem} \label{matrix-coeff}
Let $T$ be a finite lattice and let $P,Q\in\Pol_T$.
\begin{enumerate}
\item Let $\chi,\theta \in \overline{\Sur_\CL(T,P)}$ and $\tau\in\Aut(P)$.
Let also $\pi,\kappa \in \overline{\Sur_\CL(T,Q)}$ and $\sigma\in\Aut(Q)$. Then
$$f_{\chi,\tau,\theta} \, f_{\pi,\sigma,\kappa}=\begin{cases}{f_{\chi,\tau\sigma,\kappa}} & { \text{ if } \, P=Q \; \text{ and } \; \theta=\pi \,,} \\
{0} & { \text{ otherwise } \,.}\end{cases}$$
\item When $P$ varies in $\Pol_T$ and $\pi$ varies in~$\overline{\Sur_\CL(T,P)}$,
the idempotents $f_{\pi,\Id_P,\pi}$ are pairwise orthogonal.
\end{enumerate}
\fresult

\pf 
Let $\chi'=\tau^{-1}\chi\in \Sur_\CL(T,P)$, so that $j^{\chi'}=j^\chi \tau$, by Remark~\ref{Aut-action}.
If $P\neq Q$, there is no isomorphism between $P$ and $Q$, by our choice of $\Pol_T$.
Therefore we obtain $j^{\chi'}\theta j^\pi=0$, by Proposition~\ref{unique-iso}.
It follows that
$$f_{\chi,\tau,\theta} \, f_{\pi,\sigma,\kappa}=j^\chi\tau\theta j^\pi\sigma\kappa=j^{\chi'}\theta j^\pi\sigma\kappa=0 \mpoint$$

So we now assume that $P=Q$. Suppose that $f_{\chi,\tau,\theta} \, f_{\pi,\sigma,\kappa}\neq0$.
In particular, $j^\chi\tau\theta j^\pi\neq 0$, that is, $j^{\chi'}\theta j^\pi\neq 0$.
By Proposition~\ref{unique-iso}, there is a unique isomorphism $\rho:P\to P$ such that $\theta j^\pi=\rho \pmod{\Ker(j^{\chi'})}$.
Let $\theta'=\rho^{-1}\theta$ and $\chi''=\rho^{-1}\chi'$, hence $j^{\chi''}=j^{\chi'}\rho=j^\chi\tau\rho$.
Then we obtain
$$0\neq j^{\chi'}\theta j^\pi=j^\chi\tau\theta j^\pi=j^\chi\tau\rho\rho^{-1}\theta j^\pi=j^{\chi''}\theta' j^\pi \mpoint$$
Moreover, since $\theta j^\pi=\rho + h$ with $j^{\chi'}h=0$, we have $j^{\chi''}\rho^{-1}h=0$.
Therefore
$$\theta' j^\pi=\rho^{-1}\theta j^\pi=\Id_P +\rho^{-1} h =\Id_P \pmod{\Ker(j^{\chi''})} \mpoint$$
The uniqueness of the automorphism in Proposition~\ref{unique-iso} also implies that we have $\theta'(b^\pi_p)=p$ for all $p\in P$
(where $b_p^\pi=\pi\op(p)$, as before).\par

We are now in the assumptions of Proposition~\ref{compute-composite} for $\pi$, $\theta'$, and~$\chi''$.
We deduce that $\theta'=\pi$,
so that $\theta$ and $\pi$ belong to the same orbit under the action of~$\Aut(P)$.
But $\theta$ and $\pi$ belong to a chosen system of representatives $\overline{\Sur_\CL(T,P)}$.
Thus we must have $\rho=\Id_P$ and $\theta=\pi$.\par

It now follows that we can write $\theta j^\pi=\Id_P+h$, where $j^{\chi'}h=0$, that is, $j^\chi\tau h=0$. Therefore
$$f_{\chi,\tau,\theta} \, f_{\pi,\sigma,\kappa}=j^\chi\tau\theta j^\pi\sigma\kappa
=j^\chi\tau(\Id_P+h)\sigma\kappa=j^\chi\tau\sigma\kappa=f_{\chi,\tau\sigma,\kappa} \mvirg$$
as was to be shown.

(b) This follows from (a).
\endpf

%%%%%%%% Section

\section{Subalgebras corresponding to pole lattices} \label{Section-subalgebras}

\bigskip
\noindent
In this section, we show how the results of Section~\ref{Section-idempotents} imply some precise information
about the structure of the endomorphism algebra $\End_{k\CL}(T)$ of a finite lattice~$T$, where $k$ is a commutative ring.
This generalizes the results obtained in~\cite{BT3} for the special case of totally ordered lattices.

We continue to use Notation~\ref{representatives}, so $P$ is a pole lattice running through the set $\Pol_T$
and $\overline{\Sur_\CL(T,P)}$ denotes a set of representatives of $\Aut(P)$-orbits in~$\Sur_\CL(T,P)$.
Let $M_{n(T,P)}(k\Aut(P))$ denote the matrix algebra of size~$n(T,P)=|\overline{\Sur_\CL(T,P)}|$,
with rows and columns indexed by the set~$\overline{\Sur_\CL(T,P)}$, and coefficients in the group algebra $k\Aut(P)$.
If $\chi,\theta\in \overline{\Sur_\CL(T,P)}$ and $\tau\in\Aut(P)$, we let $m_{\chi,\tau,\theta}$ denote
the elementary matrix having coefficient $\tau$ in position~$(\chi,\theta)$ and zero elsewhere.

\result{Theorem} \label{matrix-algebra}
Let $T$ be a finite lattice.
For each $P\in\Pol_T$, let $\overline{\Sur_\CL(T,P)}$ be a set of representatives of the orbits for the action of the group $\Aut(P)$ on~$\Sur_\CL(T,P)$
and let $n(T,P)=|\overline{\Sur_\CL(T,P)}|$.
\begin{enumerate}
\item The map
$$\CI_T: \; \bigoplus_{P\in\Pol_T} M_{n(T,P)}(k\Aut(P)) \longrightarrow \End_{k\CL}(T) \,, \qquad m_{\chi,\tau,\theta} \mapsto f_{\chi,\tau,\theta}$$
is an algebra homomorphism (without unit elements).
\item $\CI_T$ is injective.
\item The image of $\CI_T$ is equal to the subalgebra~$\CE_T$ (without unit element) of $\End_{k\CL}(T)$ having a $k$-basis consisting of all join-morphisms $T\to T$ whose image is a pole lattice.
\end{enumerate}
\fresult

\pf
(a) Let $P,Q\in\Pol_T$. If $P\neq Q$, then $m_{\chi,\tau,\theta}$ and $m_{\pi,\sigma,\kappa}$ are not in the same block, so their product is~0,
while the product $f_{\chi,\tau,\theta}f_{\pi,\sigma,\kappa}$ is also zero.
If $P=Q$, then the relations of Theorem~\ref{matrix-coeff} are the standard relations within a matrix algebra of size~$n(T,P)$ with coefficients in the group algebra $k\Aut(P)$.\mpn

(b) Since the elements $m_{\chi,\tau,\theta}$ form a $k$-basis of $\bigoplus_{P\in\Pol_T} M_{n(T,P)}(k\Aut(P))$,
it suffices to prove that their images $f_{\chi,\tau,\theta}$ are $k$-linearly independent. Suppose that
$$\sum_{\chi,\tau,\theta} \lambda_{\chi,\tau,\theta} \, f_{\chi,\tau,\theta} =0 \mvirg$$
where $\lambda_{\chi,\tau,\theta} \in k$.
Multiplying on the left by $f_{\pi,\Id_P,\kappa}$ and on the right by $f_{\pi,\sigma,\pi}$, we are left with the terms for which $\chi=\kappa$ and $\theta=\pi$.
Therefore we obtain
$$\sum_\tau \lambda_{\kappa,\tau,\pi} \, f_{\pi,\Id_P,\kappa} \, f_{\kappa,\tau,\pi} \, f_{\pi,\sigma,\pi}
=\sum_\tau \lambda_{\kappa,\tau,\pi} \, f_{\pi,\tau\sigma,\pi}=0 \mpoint$$
Now, by Definition~\ref{def-j}, $f_{\pi,\tau\sigma,\pi}=j^\pi\tau\sigma\pi$ is a linear combination of distinct maps $T\to T$,
one of them being $j_{B^-}^\pi\tau\sigma\pi$, appearing with coefficient $\pm1$, where we use Notation~\ref{notation-B} and set $B^-=(b_e^-)_{e\in E}$.
We claim that the functions $j_{B^-}^\pi\tau\sigma\pi$ are pairwise distinct when $\tau$ varies.
This implies that each coefficient $\lambda_{\kappa,\tau,\pi}$ must be zero, proving the required linear independence.\par

To prove the claim, we write for simplicity $\rho=\tau\sigma$ and we allow $\rho$ to vary.
The group $\Aut(P)$ is isomorphic to $C_2\times C_2\times \ldots \times C_2$, where each $C_2$ acts by exchanging two twin elements of~$E_2$
and fixing the others (where $E=\Irr(P)$, as before).
So we consider some $e\in E_2$ and we let $\breve e$ be its twin. Then we get
$$j_{B^-}^\pi\rho\pi(b_e)=j_{B^-}^\pi\rho(e)=b_{\rho(e)}^-
=\begin{cases}{b_e^-} & { \text{ if } \, \rho(e)=e \,,} \\
{b_{\breve e}^-} & { \text{ if } \, \rho(e)=\breve e \,.}\end{cases}$$
We see that the functions $j_{B^-}^\pi\rho\pi$ are pairwise distinct when $\rho$ varies, proving the claim.

(c) It is clear that $\CE_T$ is a subalgebra.
Moreover, every map $j_A^\chi$ is a join-morphism, where $A$ is a family as in Lemma~\ref{order-preserving}.
Therefore $j_A^\chi\tau\theta$ is a join-morphism whose image is a pole lattice, by construction.
It follows that $f_{\chi,\tau,\theta}=j^\chi\tau\theta$ belongs to~$\CE_T$ and hence $\Im(\CI_T)\subseteq\CE_T$.\par

Now $\CE_T$ has a $k$-basis consisting of all morphisms $\varphi_{\lambda,\tau,\pi}$ described as follows.
First we fix $P\in\Pol_T$ and we let
$$\varphi_{\lambda,\tau,\pi}=\lambda\tau\pi \mvirg$$
where $\pi\in\overline{\Sur_\CL(T,P)}$, $\tau\in\Aut(P)$, $\lambda\in\overline{\Inj_\CL(P,T)}$,
and where $\overline{\Inj_\CL(P,T)}$ denotes a set of representatives of $\Aut(P)$-orbits in~$\Inj_\CL(P,T)$.
Then
$$\{\varphi_{\lambda,\tau,\pi} \mid \pi\in \overline{\Sur_\CL(T,P)}\,, \, \lambda\in \overline{\Inj_\CL(P,T)} \,, \,\tau\in\Aut(P) \}$$
is a $k$-basis of the submodule $\CE_{T,P}$ generated by all endomorphisms whose image is isomorphic to~$P$.
Allowing $P$ to vary in~$\Pol_T$, we deduce that
$$\CB=\bigcup_{P\in\Pol_T} \{\varphi_{\lambda,\tau,\pi} \mid \pi\in \overline{\Sur_\CL(T,P)}\,, \, \lambda\in \overline{\Inj_\CL(P,T)} \,, \,\tau\in\Aut(P) \}$$
is a $k$-basis of~$\CE_T=\bigoplus_{P\in\Pol_T} \CE_{T,P}$.\par

On the other hand, it follows from (a) and (b) that
$$\CB'=\bigcup_{P\in\Pol_T} \{f_{\chi,\tau,\theta} \mid \chi,\theta\in \overline{\Sur_\CL(T,P)} \,, \tau\in\Aut(P) \}$$
is a $k$-basis of~$\Im(\CI_T)$.
By Lemma~\ref{bijection-bases}, there is a bijection between $\Inj_\CL(P,T)$ and $\Sur_\CL(T,P)$.
We can also choose representatives to obtain a bijection between $\overline{\Inj_\CL(P,T)}$ and $\overline{\Sur_\CL(T,P)}$,
because $\Aut(P)$ acts freely on each side.
Therefore $\CB$ and $\CB'$ have the same cardinality.
In other words $\Im(\CI_T)$ and $\CE_T$ are free $k$-modules of the same rank.
We want to prove that the inclusion $\Im(\CI_T) \subseteq \CE_T$ is an equality (which is obvious if $k$ is a field since the dimensions are equal).\par

We now allow the base ring $k$ to vary and we write a superscript~$(k)$ to emphasize the dependence on~$k$.
Thus we have an injective algebra homomorphism
$$\CI_T^{(k)}: \; \bigoplus_{P\in\Pol_T} M_{n(T,P)}(k\Aut(P)) \longrightarrow \CE_T^{(k)}
\subseteq \End_{k\CL}(T)$$
and we let $\CX^{(k)}:=\CE_T^{(k)} / \Im(\CI_T^{(k)})$, so that we have a short exact sequence
$$\xymatrix{
0 \ar[r] & \Im(\CI_T^{(k)}) \ar[r]^-{j_k} & \CE_T^{(k)} \ar[r]^-{p_k} & \CX^{(k)} \ar[r] & 0 \mvirg
}$$
where $j_k$ is the inclusion map and $p_k$ the canonical surjection.
In the case of the ring of integers $\Z$, we see that $\CX^{(\Z)}$ is a finite abelian group,
because $\Im(\CI_T^{(\Z)})$ and $\CE_T^{(\Z)}$ are free $\Z$-modules of the same rank.
Tensoring with $k$ is right exact, so we obtain an exact sequence
$$\xymatrix{
k \otimes \Im(\CI_T^{(\Z)}) \ar[r]^-{1\otimes j_\Z} & k \otimes \CE_T^{(\Z)} \ar[r]^-{1\otimes p_\Z} & k \otimes \CX^{(\Z)}\ar[r] & 0 \mpoint
}$$
Using the canonical bases $\CB$ and~$\CB'$ of $\Im(\CI_T^{(k)})$ and $\CE_T^{(k)}$ respectively, we see that
$$k\otimes \Im(\CI_T^{(\Z)}) \cong \Im(\CI_T^{(k)}) \qquad \text{and} \qquad
k\otimes \CE_T^{(\Z)} \cong \CE_T^{(k)} \mpoint$$
Moreover the map $1\otimes j_\Z$ corresponds, under these isomorphisms, to the inclusion map~$j_k$.
In particular, considering the prime field $\F_p$ for any prime number~$p$, we obtain an exact sequence
$$\xymatrix{
\Im(\CI_T^{(\F_p)}) \ar[r]^-{j_{\F_p}} & \CE_T^{(\F_p)} \ar[r]^-{1\otimes p_\Z} & \F_p\otimes \CX^{(\Z)} \ar[r] & 0 \mpoint
}$$
Since $\F_p$ is a field and the dimensions are equal, the inclusion map $j_{\F_p}$ is an equality.
Therefore $\F_p\otimes \CX^{(\Z)}=\{0\}$ and this holds for every prime~$p$.
Thus we must have $\CX^{(\Z)}=\{0\}$, because $\CX^{(\Z)}$ is finite,
so that the inclusion map $j_\Z:\Im(\CI_T^{(\Z)}) \to \CE_T^{(\Z)}$ is an equality.
Tensoring with~$k$, it follows that the inclusion map $j_k:\Im(\CI_T^{(k)}) \to \CE_T^{(k)}$ is an equality as well, as required.
\endpf

\begin{rem}{Remark} \label{change-basis} Let $\CB$ and $\CB'$ be the two bases of $\CE_T=\Im(\CI_T)$ described in the proof.
The change of basis from $\CB$ to $\CB'$ is not obvious.
By construction, every map $j_A^\chi\tau\theta$ belongs to~$\CB$,
but beware of the fact that if $\theta$ and $\chi$ belong to $\Sur_\CL(T,P)$,
then $j_A^\chi\tau\theta$ may be a composite $T\to P' \to T$ for some pole lattice $P'$ smaller than~$P$.
This is because, in the construction of~$j_A^\chi$, the family $A= (a_e)_{e\in E}$ does not necessarily consist of distinct elements (where $E=\Irr(P)$ as before).
\end{rem}
\mpn

The image under~$\CI_T$ of the identity element of $\bigoplus_{P\in\Pol_T} M_{n(T,P)}(k\Aut(P))$ is an idempotent~$e_T$ of~$\End_{k\CL}(T)$
and $e_T$ is an identity element of~$\CE_T$. We now prove more.

\result{Theorem} \label{central-e}
For every finite lattice~$T$, let $\CE_T=\Im(\CI_T)$ be the subalgebra of~$\End_{k\CL}(T)$ appearing in Theorem~\ref{matrix-algebra},
and let $e_T$ be the identity element of~$\CE_T$.
\begin{enumerate}
\item $e_T=\displaystyle\sum_{P\in\Pol_T} \sum_{\pi\in\overline{\Sur_\CL(T,P)}} f_{\pi,\Id_P,\pi}$.
\item For any finite lattice~$T'$ and any morphism $\alpha\in\Hom_{k\CL}(T,T')$, we have $\alpha e_T = e_{T'} \alpha$.
In other words, the family of idempotents~$e_T$, for $T\in\CL$, is a natural transformation of the identity functor~$\Id_{k\CL}$.
\item $e_T$ is a central idempotent of~$\End_{k\CL}(T)$.
\item The subalgebra $\CE_T$ is a direct product factor of~$\End_{k\CL}(T)$, that is, there exists a subalgebra $\CD$ such that $\End_{k\CL}(T)=\CE_T\times\CD$ (where $\CE_T$ is identified with $\CE_T\times\{0\}$ and $\CD$ with $\{0\}\times\CD$, as usual).
\end{enumerate}
\fresult

\pf (a) The identity element of $\displaystyle\bigoplus_{P\in\Pol_T} M_{n(T,P)}(k\Aut(P))$ is equal to
$$\sum_{P\in\Pol_T} \sum_{\pi\in\overline{\Sur_\CL(T,P)}} m_{\pi,\Id_P,\pi} \mpoint$$
Taking its image under~$\CI_T$ yields the required formula.\mpn

(b) We have seen in the proof of Theorem~\ref{matrix-algebra} that every element of the canonical basis~$\CB$ of~$\CE_T$
has the form $\varphi_{\lambda,\tau,\pi}=\lambda\tau\pi$,
where $\pi\in\overline{\Sur_\CL(T,P)}$, $\lambda\in\overline{\Inj_\CL(P,T)}$ and $\tau\in\Aut(P)$.
Passing to the opposite, we obtain
$$\varphi_{\lambda,\tau,\pi}\op=\pi\op\tau\op\lambda\op=\pi\op\tau^{-1}\lambda\op$$
with $\pi\op \in \Inj_\CL(P\op,T\op)$ and $\lambda\op\in \Sur_\CL(T\op,P\op)$.
It follows that the opposite of the canonical basis element $\varphi_{\lambda,\tau,\pi}$ of~$\CE_T$
is the canonical basis element $\varphi_{\pi\op,\tau^{-1},\lambda\op}$ of~$\CE_{T\op}$.
Therefore, the opposite of the identity element~$e_T$ of~$\CE_T$ must belong to~$\CE_{T\op}$.
Moreover, it must be the identity element of~$\CE_{T\op}$, because taking opposites behaves well with respect to composition.
Therefore $(e_T)\op=e_{T\op}$.\par

Now if $\alpha:T\to T'$ is a join-morphism (i.e. $\alpha$ is in~$\CL$), then the image of a pole sublattice of~$T$ is a pole sublattice of~$T'$.
It follows that composition with~$\alpha$ maps $e_T$ to a linear combination of join-morphisms with a pole lattice as an image,
hence invariant under the idempotent element~$e_{T'}$. In other words, we have
$$\alpha \,e_T=e_{T'} \,\alpha \,e_T \mpoint$$
Applying this equation to ${T'}\op$, $T\op$, and the morphism $\alpha\op:{T'}\op \to T\op$, we obtain $\alpha\op \,e_{{T'}\op}=e_{T\op} \,\alpha\op \,e_{{T'}\op}$.
Passing to opposites and using the above equality $(e_T)\op=e_{T\op}$, we get
$$e_{T'} \, \alpha \, e_T = e_{T'} \, \alpha \mpoint$$
The two displayed equations yield $\alpha e_T = e_{T'} \alpha$.
This holds as well if $\alpha$ is replaced by a $k$-linear combination of join-morphisms (i.e. $\alpha$ is in~$k\CL$), as was to be shown.\mpn

(c) This is a special case of~(b).\mpn

(d) This follows immediately from~(c).
\endpf

%%%%%%%%%%%%%%%%%

\section{Correspondence functors for pole lattices} \label{Section-functors-pole}

\bigskip
\noindent
In this section, we first consider the special case of the endomorphism algebra of a pole lattice~$Q$.
We determine completely the structure of this algebra.
Applying the fully faithful functor $T\mapsto F_T$, we deduce a direct sum decomposition of the correspondence functor~$F_Q$,
providing an explicit description of~$F_Q$ for any pole lattice~$Q$.
In particular, when $k$ is a field of characteristic different from~2, $F_Q$ is semi-simple.
At the end the section, we return to an arbitrary finite lattice~$T$ and describe direct summands of~$F_T$ corresponding to pole lattices inside~$T$.
The results are generalizations of those obtained in~\cite{BT3} in the special case of totally ordered lattices.

\result{Theorem} \label{endom-pole} Let $Q$ be a pole lattice.
\begin{enumerate}
\item The homomorphism of $k$-algebras of Theorem~\ref{matrix-algebra}
$$\CI_{Q}:\bigoplus_{P\in\Pol_Q} M_{n(Q,P)}(k\Aut(P))\longrightarrow \End_{k\CL}(Q) \,, \qquad
m_{\chi,\tau,\theta}\mapsto f_{\chi,\tau,\theta} \mvirg$$
is an isomorphism.
\item In particular, if $k$ is a field and if either $Q$ is totally ordered or if $k$ is a field of characteristic different from~2,
then $\End_{k\CL}(Q)$ is semi-simple.
\end{enumerate}
\fresult

\pf
(a) Since any join-morphism $\varphi:Q\to Q$ has an image which is a pole lattice,
the subalgebra $\CE_Q$ of~$\End_{k\CL}(Q)$ appearing in Theorem~\ref{matrix-algebra} is the whole of~$\End_{k\CL}(Q)$.
Therefore, the homomorphism $\CI_{Q}$ is surjective.
By Theorem~\ref{matrix-algebra}, $\CI_{Q}$ is injective, hence an isomorphism.\mpn

(b) If $Q$ is totally ordered, then so is each~$P$ and $\Aut(P)$ is the trivial group.
Thus we get matrix algebras $M_{n(Q,P)}(k)$. 
If $Q$ is not totally ordered, then each $\Aut(P)$ is a 2-group (and at least one of them is nontrivial, namely $\Aut(Q)$).
The group algebra $k\Aut(P)$ is semi-simple when the characteristic of~$k$ is different from~2 (Maschke's theorem).
Therefore any matrix algebra $M_q(k\Aut(P))$ is semi-simple and it follows that the direct sum is semi-simple as well.
\endpf

Now we consider the central idempotents of $\End_{k\CL}(Q)$ corresponding to the above decomposition into matrix algebras.

\result{Notation} \label{beta}
For any pole lattice~$P\in\Pol_Q$, set 
$$\beta_{Q,P}:=\sum_{\pi\in\overline{\Sur_\CL(Q,P)}} f_{\pi,\Id_Q,\pi}\mpoint$$
In particular, when $P=Q$, then $\Sur_\CL(Q,Q)=\Aut(Q)$ and $\overline{\Sur_\CL(Q,Q)}$ is a singleton which can be chosen to be $\{\Id_{Q}\}$.
We then define
$$\varepsilon_Q:=\beta_{Q,Q}=f_{\Id_Q,\Id_Q,\Id_Q}=j^{\Id_Q}= \sum_{\emptyset \subseteq Y \subseteq E} (-1)^{|E-Y|} \rho_Y\mvirg$$
using Proposition~\ref{rho_Y}, with $E=\Irr(Q)$ and $\rho_Y \in\End_\CL(Q)$ defined by~(\ref{def-rho-Y}).
\fresult

\result{Proposition} \label{central-idempotents} 
The elements~$\beta_{Q,P}$, for $P\in\Pol_Q$, are orthogonal central idempotents of $\End_{k\CL}(Q)$, and their sum is equal to the identity.
In particular, the central idempotent $\varepsilon_Q$ satisfies
$$\varepsilon_Q\End_{k\CL}(Q)\cong k\Aut(Q) \mpoint$$
\fresult

\pf
For every $\pi\in\overline{\Sur_\CL(Q,P)}$, the inverse image of $f_{\pi,\Id_P,\pi}$ under the algebra isomorphism~$\CI_Q$
of Theorem~\ref{endom-pole} is the matrix $m_{\pi,\Id_P,\pi}$ of the component $M_{n(Q,P)}(k\Aut(P))$ indexed by~$P$. 
Summing over all $\pi\in\overline{\Sur_\CL(Q,P)}$, it follows that the inverse image of~$\beta_{Q,P}$ under~$\CI_Q$
is the identity element of $M_{n(Q,P)}(k\Aut(P))$.
The first statement follows.\par

In the case $P=Q$, we know that $\overline{\Sur_\CL(Q,Q)}$ is a singleton, so that the corresponding matrix algebra has size~1.
The inverse image of $\varepsilon_Q$ under~$\CI_Q$
is the identity element $m_{\Id_Q,\Id_Q,\Id_Q}$ of the component $M_1(k\Aut(Q))\cong k\Aut(Q)$.
Clearly $\varepsilon_Q\End_{k\CL}(Q)\cong M_1(k\Aut(Q))\cong k\Aut(Q)$.
\endpf

We want to use the fully-faithful functor $F_?: k\CL \to \CF_k$ (see Theorem~\ref{fully-faithful})
to deduce information on the correspondence functor~$F_Q$.
We already know that $F_Q$ is projective, because the pole lattice~$Q$ is distributive (see Theorem~\ref{fully-faithful}).
We apply the functor $F_?: k\CL \to \CF_k$ to the map $j^\pi\in \Hom_{k\CL}(P,Q)$ defined in~(\ref{def-j}), where $\pi\in\overline{\Sur_\CL(Q,P)}$.
By Proposition~\ref{vanish-on-H} we obtain a morphism
$$F_{j^\pi}:F_{P} \longrightarrow F_{Q}$$
which vanishes on~$H_P$, where $H_P$ is defined by~(\ref{def-H}).
By Proposition~\ref{ker-j}, this induces an injective morphism
$$\overline F_{j^\pi}:F_{P}/H_{P} \longrightarrow F_{Q}$$
which embeds $F_{P}/H_{P}$ as a direct summand of~$F_{Q}$, corresponding to the idempotent $f_{\pi,\Id_P,\pi}=j^\pi\pi$.
In particular, for $P=Q$, we have $f_{\Id_Q,\Id_Q,\Id_Q}=j^{\Id_Q}=\varepsilon_Q$ and we obtain an idempotent endomorphism
$F_{\varepsilon_Q}$ of~$F_{Q}$ with kernel~$H_{Q}$.

\result{Theorem} \label{functor-SQ}
Let $Q$ be a pole lattice and define $\S_Q:=F_{Q}/H_{Q}$, where $H_Q$ is defined by~(\ref{def-H}).
\begin{enumerate}
\item $\S_Q$ is a projective correspondence functor.
\item There are isomorphisms of correspondence functors
\begin{eqnarray*}
F_{\varepsilon_Q}F_Q&\cong& \S_Q \mvirg \\
F_{\beta_{Q,P}}F_Q&\cong& \S_P^{n(Q,P)} \mvirg\; \text{ for each }\; P\in \Pol_Q\mvirg \\
F_Q&\cong&\bigoplus_{P\in\Pol_Q} \S_P^{n(Q,P)}\mpoint
\end{eqnarray*}
\end{enumerate}
\fresult

\pf (a) Since the pole lattice $Q$ is distributive, $F_Q$ is projective (Theorem~\ref{fully-faithful}).
Therefore so is its direct summand $\S_Q$.\mpn

(b) Since the functor $F_?: k\CL \to \CF_k$ is fully faithful (Theorem~\ref{fully-faithful}), it induces an isomorphism of $k$-algebras 
$$\End_{k\CL}(Q)\cong\End_{\CF_k}(F_{Q})\mpoint$$
Now the idempotents $f_{\pi,\Id_P,\pi}$ of $\End_{k\CL}(Q)$, for $\pi\in\overline{\Sur_\CL(Q,P)}$ and $P\in\Pol_Q$,
are orthogonal and their sum is equal to the identity, by Theorem~\ref{endom-pole}.
It follows that the endomorphisms $F_{f_{\pi,\Id_P,\pi}}$ of~$F_{Q}$ are orthogonal idempotents, and their sum is the identity.
Hence we obtain a decomposition of correspondence functors
$$F_Q=\bigoplus_{\substack{ P\in\Pol_Q \\ \pi\in\overline{\Sur_\CL(Q,P)} }} F_{f_{\pi,\Id_P,\pi}}\big(F_Q\big)\mpoint$$
By surjectivity of $\pi:Q \to P$, the image of $F_{f_{\pi,\Id_P,\pi}}=F_{j^\pi}F_{\pi}:F_Q \to F_Q$ is equal to the image of~$F_{j^\pi}:F_P \to F_Q$.
Therefore $F_{f_{\pi,\Id_P,\pi}}\big(F_Q\big)=F_{j^\pi}\big(F_P\big)$.
By Proposition~\ref{ker-j}, the image $F_{j^\pi}\big(F_P\big)$ is isomorphic to~$\S_P=F_P/H_P$ and it follows that
$$F_{f_{\pi,\Id_P,\pi}}\big(F_Q\big) \cong \S_P \mpoint$$
Taking $P=Q$ and $f_{\Id_Q,\Id_Q,\Id_Q}=j^{\Id_Q}=\varepsilon_Q$, we obtain the first isomorphism $F_{\varepsilon_Q}F_Q\cong \S_Q$.
Summing over all $\pi\in\overline{\Sur_\CL(Q,P)}$ for a fixed~$P$, we obtain the second isomorphism.
Finally, summing over all $P\in\Pol_Q$ and all $\pi\in\overline{\Sur_\CL(Q,P)}$, we obtain the third isomorphism.
\endpf

\result{Corollary} \label{Hom-zero} Let $P$ and $P'$ be pole lattices. Then
$$\Hom_{\CF_k}(\S_P,\S_{P'})
\cong \left\{\begin{array}{ll} \{0\} &\hbox{if}\;P\not\cong P' \mvirg\\
k\Aut(P)&\hbox{if}\; P\cong P' \mpoint \end{array}\right.$$
\fresult

\pf 
Since $\S_P\cong F_{\varepsilon_P}F_P$, the case $P=P'$ follows from Proposition~\ref{central-idempotents}.
Now if $P\not\cong P'$, it is easy to choose a large enough pole lattice~$Q$ such that $\Sur_\CL(Q,P)\neq\emptyset$ and $\Sur_\CL(Q,P')\neq\emptyset$.
Using the central idempotents $\beta_{Q,P}$ and $\beta_{Q,P'}$ of~Proposition~\ref{central-idempotents}, we obtain
$$\Hom_{\CF_k}(F_{\beta_{Q,P}}F_Q,F_{\beta_{Q,P'}}F_Q)\cong \Hom_{\CF_k}(\S_P,\S_{P'})^{n(Q,P)\cdot n(Q,P')}\mpoint$$
Since $F_{\beta_{Q,P}}$ and $F_{\beta_{Q,P'}}$ are central idempotents of $\End_{\CF_k}(F_{Q})$, and since they are orthogonal if $P\not\cong P'$,
it follows that $\Hom_{\CF_k}(F_{\beta_{Q,P}}F_Q,F_{\beta_{Q,P'}}F_Q)=0$ if $P\not\cong P'$, hence $\Hom_{\CF_k}(\S_P,\S_{P'})=\{0\}$.
\endpf

\result{Remark} {\rm
Corollary~\ref{Hom-zero} actually holds for the fundamental functors associated with any two finite posets.
This more general result will be proved in another paper.
}
\fresult

\bigskip
Now we prove that the functor $\S_Q$ is actually isomorphic to a fundamental functor and we compute the ranks of all its evaluations.

\result{Theorem} \label{S_Q-fundamental} Let $Q$ be a pole lattice and let $R$ be the order relation on the set $E=\Irr(Q)$ of irreducible elements of~$Q$.
Let $\S_Q=F_{Q}/H_{Q}$.
\begin{enumerate}
\item $\S_Q$ is isomorphic to the fundamental functor $\S_{E,R\op}$.
\item For any finite set $X$, the $k$-module $\S_Q(X)$ is free of rank 
$${\rm rank}\big(\S_Q(X)\big)=\sum_{i=0}^{|E|}(-1)^i{|E|\choose i}(|Q|-i)^{|X|} \mpoint$$
\end{enumerate}
\fresult

\pf 
(a) We use the element $\gamma_{Q\op} \in F_Q(E^0)$ defined in~(\ref{def-gamma}), where $E^0=\Irr(Q\op)$.
By a well-known result of lattice theory (Theorem~6.2 in~\cite{Ro}), the distributive lattice $Q\op$ is isomorphic to $\Idown(E^0,R^0)$,
where $R^0$ is the order relation on~$E^0$ viewed as a subset of~$Q\op$,
so that $(E^0,R^0)$ is the poset of irreducible elements in~$Q\op$.
Note that the isomorphism $Q\op \cong \Idown(E^0,R^0)$ can also be checked easily and directly because $Q\op$ is a pole lattice.
Recall that
$$\gamma_{Q\op}=\sum_{A\subseteq E^0} (-1)^{|A|} \eta_A^\circ \mvirg$$
where $\eta_A^\circ:E^0\to Q$ denotes the same map as $\eta:E^0\to Q\op$ and $\eta$ is defined by
$$\forall e^0\in E^0,\;\;\eta_A(e^0)=\left\{\begin{array}{ll} s(e^0)&\hbox{if}\;e^0\in A\mvirg \\
e^0&\hbox{if}\;e^0\notin A \mvirg \end{array}\right.
$$
because $r(e^0)$ in the lattice $Q\op$ is equal to~$s(e^0)$ in the lattice~$Q$.\par

Now we define $\omega:E^0\to Q$ by
$$\omega(e^0)=\left\{\begin{array}{ll} s(e^0)&\hbox{if}\;e^0\in E_1^0\mvirg \\
e^0&\hbox{if}\;e^0\in E_2^0 \mvirg
\end{array}\right.$$
and we notice that $\omega$ is actually a bijection between $E^0$ and~$E=\Irr(Q)$,
because in a pole lattice we have $E_1=s(E_1^0)$ and $E_2=E_2^0$ (by an easy application of Lemma~\ref{irred-pole}).
Then $\omega\in F_{Q}(E^0)$ and when we apply the idempotent $F_{\varepsilon_Q}$ we claim that we obtain
\begin{equation}\label{omega-gamma}
F_{\varepsilon_Q} (\omega)=(-1)^{|E_1|} \gamma_{Q\op} \mpoint
\end{equation}
The definition of $\varepsilon_Q$ (see Notation~\ref{beta}) yields
$$F_{\varepsilon_Q} (\omega)=\varepsilon_Q\,\omega= \sum_{Y\subseteq E} (-1)^{|E-Y|} \rho_Y \, \omega \mpoint$$
The definition of $\rho_Y$ in~(\ref{def-rho-Y}) splits into two cases. If $e^0\in E_1^0$, then
$$(\rho_Y \, \omega)(e^0)=\rho_Y(s(e^0))
=\left\{\begin{array}{ll} s(e^0)&\hbox{if}\; s(e^0)\in Y \mvirg \\
r(s(e^0))=e^0&\hbox{if}\;s(e^0)\notin Y \mpoint \end{array}\right.$$
If now $e^0\in E_2^0$, then
$$(\rho_Y \, \omega)(e^0)=\rho_Y(e^0)
=\left\{\begin{array}{ll} e^0&\hbox{if}\; e^0\in Y \mvirg \\
s(e^0)&\hbox{if}\;e^0\notin Y \mpoint \end{array}\right.$$
For each $Y\subseteq E$, we define $A\subseteq E^0$ by
$$Y\cap E_1=s(A\cap E_1^0) \qquad \text{ and } \qquad Y\cap E_2=E_2-(A\cap E_2) \mpoint$$
Thus we have decompositions
$$Y=(Y\cap E_1)\sqcup (Y\cap E_2)\subseteq E  \qquad \text{ and } \qquad A =(A\cap E_1^0) \sqcup (A\cap E_2^0) \subseteq E^0$$
and $A$ runs through all subsets of~$E^0$ when $Y$ runs through all subsets of~$E$.
If $e^0\in E_1^0$, then
$$s(e^0) \in Y\cap E_1 \iff e^0\in A\cap E_1^0$$
while if $e^0\in E_2^0$, then
$$e^0 \notin Y\cap E_2 \iff e^0\in A\cap E_2^0 \mpoint$$
Therefore the two cases merge into one and we obtain
$$(\rho_Y \, \omega)(e^0)=\left\{\begin{array}{ll} e^0&\hbox{if}\; e^0\notin A \mvirg \\
s(e^0)&\hbox{if}\;e^0\in A \mvirg \end{array}\right.$$
so that $\rho_Y \, \omega = \eta_A^\circ$.\par

As far as the signs are concerned, we have
$$|E-Y| = |E_1-(Y\cap E_1)| + |E_2-(Y\cap E_2)|
= |E_1^0-(A\cap E_1^0)| + |A\cap E_2^0| \mvirg$$
hence
$$(-1)^{|E-Y|} = (-1)^{|E_1^0|} \cdot (-1)^{|A\cap E_1^0|} \cdot (-1)^{|A\cap E_2^0|} = (-1)^{|E_1|} \cdot (-1)^{|A|} \mpoint$$
It now follows that
$$\begin{array}{rcl}
F_{\varepsilon_Q} (\omega)&=&
\displaystyle\sum_{Y\subseteq E} (-1)^{|E-Y|} \rho_Y \, \omega \\
&=&(-1)^{|E_1|}\displaystyle\sum_{A\subseteq E^0} (-1)^{|A|} \eta_A^\circ \\
&=&(-1)^{|E_1|}\gamma_{Q\op} \mpoint
\end{array}$$
This proves Claim~\ref{omega-gamma} above.\par

Now $F_{Q}$ is generated by $\omega\in F_{Q}(E^0)$, because it is generated by $\iota\in F_{Q}(E)$
(where $\iota:E\to Q$ is the inclusion), hence also by any injection from the set~$E^0$ to~$Q$,
by composing~$\iota$ with a bijection between $E^0$ and~$E$.
Since $F_{\varepsilon_Q}$ is an idempotent endomorphism of the correspondence functor $F_{Q}$,
we see that $F_{\varepsilon_Q}F_{Q}$ is generated by $F_{\varepsilon_Q} (\omega)$.
In other words, in view of Claim~\ref{omega-gamma} above, $F_{\varepsilon_Q}F_{Q}$ is generated by~$\gamma_{Q\op}
\in F_{Q}(E^0)$.
Now Theorem~\ref{SERgamma} asserts that the subfunctor of~$F_Q$ generated by
$\gamma_{Q\op}$ is isomorphic to~$\S_{E^0,R^0}$, where $(E^0,R^0)$ is the poset of irreducible elements in~$Q\op$.
But $(E^0,R^0)\cong (E,R\op)$ via the map $\omega:E^0\to E$ described above.
Therefore, using the isomorphism of Theorem~\ref{functor-SQ},  we obtain
$$\S_Q\cong F_{\varepsilon_Q}F_{Q}=\langle \gamma_{Q\op}\rangle
\cong \S_{E^0,R^0} \cong \S_{E,R\op} \mpoint$$

(b) By Definition~\ref{def-H}, the canonical $k$-basis of $\S_Q(X)=F_Q(X)/H_Q(X)$ is the set $Z(X)$ of all maps $\varphi:X\to Q$
such that $E\subseteq \varphi(X) \subseteq Q$.
Therefore $\S_Q(X)$ is free of rank $|Z(X)|$.
The number of maps in~$Z(X)$ has been computed in Lemma~8.1 of~\cite{BT2} and the formula is actually well-known.
The formula asserts that this rank is equal to
$$|Z(X)|=\sum_{i=0}^{|E|}(-1)^i{|E|\choose i}(|Q|-i)^{|X|}$$
as required.
\endpf

\begin{rem}{Remark} \label{proj-fundamental-2}
In view of the projectivity of $\S_Q$ (Theorem~\ref{functor-SQ}),
the isomorphism $\S_Q\cong \S_{E,R\op}$ provides another proof of the projectivity of the fundamental functor $\S_{E,R\op}$
whenever $(E,R\op)$ is a pole poset.
This was first proved in Theorem~\ref{proj-fundamental}.
\end{rem}

\begin{rem}{Remark}
The formula for the rank in Theorem~\ref{S_Q-fundamental} is a special case of the general formula
proved in~\cite{BT4} for the rank of the evaluation of any fundamental functor.
We have given here a direct proof in the case of a pole lattice because it is easy, while the proof in the general case is much more elaborate.
\end{rem}

\bigskip
When $k$ is a field, we get even more.

\result{Corollary} \label{simple-SQ}
Let $k$ be a field. Let $Q$ be a pole lattice and let $(E,R)$ be the poset of irreducible elements in~$Q$.
\begin{enumerate}
\item For any simple $k\Aut(Q)$-module~$V$, the functor $\S_Q\otimes_{k\Aut(Q)}V$ is simple, isomorphic to $S_{E,R\op,V}$.
\item The correspondence functor $\S_Q$ is projective and injective.
\item If either $\Aut(Q)$ is trivial (which occurs if $Q$ is totally ordered) or if the characteristic of~$k$ is different from~2,
the correspondence functor $S_{E,R\op,V}$ is simple, projective, and injective.
\item Under the assumption of~(c), $\S_Q$ decomposes as a direct sum of simple (and projective) functors
$$\S_Q \cong \bigoplus_V S_{E,R\op,V} \mvirg$$
where $V$ runs over simple $k\Aut(Q)$-modules up to isomorphism.
\item Under the assumption of~(c), $F_Q$ decomposes as a direct sum of simple (and projective) functors
$$F_Q \cong \bigoplus_{P\in\Pol_Q} \bigoplus_{V_P} \, (S_{E_P,R_P\op,V_P})^{n(Q,P)} \mvirg$$
where $(E_P,R_P)$ denotes the poset of irreducible elements in~$P$ and
where $V_P$ runs over simple $k\Aut(P)$-modules up to isomorphism.

\end{enumerate}
\fresult

\pf 
(a) Using Lemma~\ref{irred-pole}, it is easy to check that $\Aut(Q)=\Aut(E,R)=\Aut(E,R\op)$, so $V$ is a $k\Aut(E,R\op)$-module.
Recall that the fundamental correspondence functor $S_{E,R\op}$ has a right $k\Aut(E,R\op)$-module structure
(in the sense that each evaluation $S_{E,R\op}(X)$ is a right $k\Aut(E,R\op)$-module, in a compatible way with all morphisms, which act on the left).
Moreover, by Proposition~\ref{precursor}, we know that the simple functor $S_{E,R\op,V}$ is obtained from the fundamental functor $S_{E,R\op}$ by simply tensoring with~$V$~:
$$\S_{E,R\op,V} \cong S_{E,R\op} \otimes_{\Aut(E,R\op)} V \,,
\qquad\text{ that is,} \qquad
\S_{E,R\op,V} \cong \S_Q\otimes_{k\Aut(Q)}V \mvirg$$
as required.\mpn

(b) $\S_Q$ is projective by Theorem~\ref{functor-SQ}. Since $k$ is a field, it is also injective by Theorem~10.6 in~\cite{BT2}.\mpn

(c) When either $\Aut(Q)$ is trivial or the characteristic of~$k$ is different from~2, $k\Aut(Q)$ is semi-simple and every simple $k\Aut(Q)$-module is projective.
Moreover, every simple $k\Aut(Q)$-module has dimension~1
because $\Aut(Q)$ is an elementary abelian 2-group (the only roots of unity needed are~$\pm1$).
Therefore we have an isomorphism of $k\Aut(Q)$-modules
$$k\Aut(Q) \cong \bigoplus_{V \text{ simple}} V \mvirg$$
where $V$ runs over all simple $k\Aut(Q)$-modules up to isomorphism.
It follows that
$$\S_Q\cong \S_Q\otimes_{k\Aut(Q)} k\Aut(Q) \cong \bigoplus_{V \text{ simple}} \S_Q\otimes_{k\Aut(Q)}V
\cong \bigoplus_{V \text{ simple}} S_{E,R\op,V} \mpoint$$
Since $\S_Q$ is projective and injective by~(b), so is each of its simple direct summands $S_{E,R\op,V}$.\mpn

(d) The decomposition of~$\S_Q$ was proved above.\mpn

(e) The decomposition of~$F_Q$ follows immediately from (d) and Theorem~\ref{functor-SQ}.
\endpf

In the special case of totally ordered lattices, the results of Corollary~\ref{simple-SQ} were already obtained in Corollary~11.11 of~\cite{BT3}.
Also, notice that (c) provides another proof of the implication (b) $\Rightarrow$ (a) in Theorem~\ref{characterization-proj}.

Our last purpose in this section is to find, for any finite lattice~$T$, all the direct summands of $F_T$ isomorphic to a functor~$\S_P$ corresponding to a pole lattice~$P$.
Recall that $e_T$ denotes the central idempotent of~$\End_{k\CL}(T)$ which is an identity element for the subalgebra $\CE_T$ (see Theorem~\ref{central-e}).

\result{Theorem} \label{pole-direct-factors}
Let $T$ be a finite lattice.
For every finite set $X$, let $F_T^{\rm pole}(X)$ be the $k$-submodule of~$F_T(X)$ generated by all the maps $\varphi:X\to T$ such that $\varphi(X)$ is a pole subposet of~$T$.
\begin{enumerate}
\item $F_T^{\rm pole}=F_{e_T}(F_T)$ and this is a subfunctor of~$F_T$.
\item $F_T^{\rm pole}$ is a projective direct summand of~$F_T$, isomorphic~to
$$F_T^{\rm pole}\cong \bigoplus_{\substack{ P\in\Pol_T \\ \pi\in\overline{\Sur_\CL(T,P)} }} \S_P
=\bigoplus_{P\in\Pol_T} \S_P^{n(T,P)}\mpoint$$
\item If $Q$ is a pole lattice, the image of a join-morphism $F_{Q} \to F_T$ is contained in~$F_T^{\rm pole}$.
In particular, any subfunctor of $F_T$ isomorphic to the functor~$\S_Q$ is contained in~$F_T^{\rm pole}$.
\item $\Hom_{\CF_k}\big(F_T^{\rm pole}, F_{\Id-e_T}(F_T)\big)=\{0\}$ and
$\Hom_{\CF_k}\big(F_{\Id-e_T}(F_T), F_T^{\rm pole}\big)=\{0\}$.
\item The splitting of the surjection $F_{e_T}:F_T\to F_T^{\rm pole}$ is natural in~$T$.
\end{enumerate}
\fresult

\pf
(a) Let $\varphi:X\to T$ be a map such that $\varphi(X)$ is a pole subposet of~$T$.
Let $Q$ be the join-closure of~$\varphi(X)$, so that $\varphi=j\psi$, where $\psi:X\to Q$ and where $j:Q\to T$ is the inclusion map.
It is easy to see that $Q$ is a join-closed pole lattice.
Thus $j\in\Hom_{k\CL}(Q,T)$ and so $j=je_Q$ because $e_Q\in \End_{k\CL}(Q)$ is the identity morphism by Theorem~\ref{endom-pole}.
Now $je_Q=e_Tj$ by Theorem~\ref{central-e}, hence $j=e_Tj$.
Therefore
$$\varphi=j\psi=e_T j\psi = e_T\varphi =F_{e_T}(\varphi) \mvirg$$
proving that $\varphi\in F_{e_T}(F_T(X))$.\par

Conversely, if $\varphi\in F_{e_T}(F_T(X))$, then we can write $\varphi= F_{e_T}(\psi)=e_T\psi$ for some map $\psi:X\to T$.
Since $e_T$ is, by construction, a linear combination of maps whose image is a pole poset, so is $e_T\psi$,
proving that $\varphi\in F_T^{\rm pole}(X)$.\par

This shows that $F_T^{\rm pole}=F_{e_T}(F_T)$ and the latter is a subfunctor of~$F_T$.\mpn

(b) As in the proof of Theorem~\ref{functor-SQ}, we apply the fully faithful functor $k\CL\to\CF_k$ defined by $T\mapsto F_T$.
There is a direct sum decomposition of functors
$$F_T=F_{e_T}(F_T) \oplus F_{\Id-e_T}(F_T) = F_T^{\rm pole} \oplus F_{\Id-e_T}(F_T) \mpoint$$
The idempotent $e_T$ is the sum of the orthogonal idempotents $f_{\pi,\Id_P,\pi}$ of~$\End_{k\CL}(T)$, for $\pi\in\overline{\Sur_\CL(T,P)}$ and $P\in\Pol_T$.
It follows that the endomorphisms $F_{f_{\pi,\Id_P,\pi}}$ of~$F_T$ are orthogonal idempotents with sum~$F_{e_T}$.
Hence we obtain a direct sum decomposition of correspondence functors
$$F_T^{\rm pole}=F_{e_T}(F_T)
=\bigoplus_{\substack{ P\in\Pol_T \\ \pi\in\overline{\Sur_\CL(T,P)} }} F_{f_{\pi,\Id_P,\pi}}\big(F_T\big) \mpoint$$
By Proposition~\ref{ker-j}, the image of~$F_{f_{\pi,\Id_P,\pi}}=F_{j^\pi\pi}$ is isomorphic to~$F_P/H_P=\S_P$
and is projective by Theorem~\ref{functor-SQ}, proving the result.\mpn

(c) Let $\alpha: F_{Q} \to F_T$ be a morphism of correspondence functors where $Q$ is a pole lattice.
Since the functor $T\mapsto F_T$ is full, $\alpha$ is the image of a morphism $Q\to T$ in~$k\CL$,
which is in turn a linear combination of join-morphisms $f:Q\to T$.
Any such~$f$ has an image which is a pole subposet of~$T$.
Therefore, for any function $\varphi:X\to Q$, the image of $f\varphi$ is a pole subposet of~$T$.
It follows that the image of the map $F_f(X):F_{Q}(X) \to F_T(X)$ is contained in $F_T^{\rm pole}(X)$.
Therefore, the image of the map $F_f$ is contained in~$F_T^{\rm pole}$
and so the image of~$\alpha$ is contained in~$F_T^{\rm pole}$.\par

The special case follows from the fact that $\S_{Q}$ is a subfunctor of~$F_{Q}$, by Theorem~\ref{functor-SQ}.\mpn

(d) The first statement is a consequence of (b) and~(c), while the second one follows from a dual argument. Details are left to the reader.\mpn

(e) By Theorem~\ref{central-e}, the family of idempotents~$e_T$, for $T\in\CL$, is a natural transformation of the identity functor~$\Id_{k\CL}$. Therefore the family of idempotents~$F_{e_T}$, for $T\in\CL$, is a natural transformation of the identity functor~$\Id_{\CF_k}$.
\endpf

\result{Corollary}
Let $F$ be a correspondence functor and let $F^{\rm pole}$ be the sum of all the images of morphisms $F_P\to F$, where $P$ varies among pole lattices.
\begin{enumerate}
\item The subfunctor $F^{\rm pole}$ is the image of an idempotent natural transformation $\varepsilon_F: F\to F$, so that $F^{\rm pole}$ is a direct summand of~$F$.
\item $\Hom_{\CF_k}\big(F^{\rm pole}, (\Id-\varepsilon_F)(F)\big)=\{0\}$.
\item The idempotent $\varepsilon_F$ is natural in~$F$. In other words, when $F$ varies among correspondence functors,
the family of idempotents $\varepsilon_F$ is a natural transformation of the identity functor $\CF_k\to\CF_k$.
\end{enumerate}
\fresult

\pf
We only sketch the main arguments of the proof.
By Yoneda's lemma applied to a set of generators of~$F$, there is some index set~$I$ and a surjective morphism from a direct sum of representable functors
$$\pi:\bigoplus_{i\in I} k\CC(-,E_i) \longrightarrow F$$
and each $k\CC(-,E_i)$ is projective.
Moreover, $k\CC(-,E_i)$ is isomorphic to~$F_{T_i}$ for some distributive lattice~$T_i$ (by Lemma~\ref{iso-complement}).
It follows that there is an exact sequence
$$\xymatrix{
\displaystyle \bigoplus_{j\in J} F_{U_j} \ar[r]& \displaystyle\bigoplus_{i\in I} F_{T_i} \ar[r]^{\quad\pi}& F \ar[r]& 0 }$$
where $U_j$ is again a distributive lattice for each $j$ in some index set~$J$.
Let us write $\varepsilon$ for the direct sum of the idempotent endomorphisms of Theorem~\ref{pole-direct-factors}, independently of the lattices involved.
Thus we have a commutative diagram
$$\xymatrix{
\displaystyle \bigoplus_{j\in J} F_{U_j}  \ar[d]^-{\varepsilon} \ar[r]&  \displaystyle\bigoplus_{i\in I} F_{T_i} \ar[d]^-{\varepsilon} \ar[r]^{\quad\pi}
&F \ar[d]^-{\varepsilon_F} \ar[r] &0 \\
\displaystyle\bigoplus_{j\in J} F_{U_j} \ar[r]& \displaystyle\bigoplus_{i\in I} F_{T_i} \ar[r]^-{\pi}& F \ar[r]&0
}$$
where $\varepsilon_F: F\to F$ is induced by~$\varepsilon$.
It is easy to check that $\varepsilon_F$ is an idempotent morphism and that $\Im(\varepsilon_F)\subseteq F^{\rm pole}$,
because $\Im(\varepsilon_F)=\Im(\pi\varepsilon)$ and this is the image under~$\pi$ of correspondence functors associated to pole lattices, by Theorem~\ref{pole-direct-factors}.
Moreover, any pole lattice~$P$ is distributive, so $F_P$ is projective.
Therefore any morphism $F_P\to F$ lifts to a morphism $F_P\to\bigoplus_{i\in I} F_{T_i}$ whose image must be contained in~$\Im(\varepsilon)$.
Thus $F^{\rm pole}$ is contained in $\pi(\Im(\varepsilon))=\Im(\varepsilon_F)$.\par

The proofs of~(b) and (c) are similar.
\endpf

%%%%%%%%%%%%%%%%%

\bigskip\bigskip
\noindent
Serge Bouc, CNRS-LAMFA, Universit\'e de Picardie - Jules Verne,\\
33, rue St Leu, F-80039 Amiens Cedex~1, France.\\
{\tt serge.bouc@u-picardie.fr}

\medskip
\noindent
Jacques Th\'evenaz, Institut de math\'ematiques, EPFL, \\
Station~8, CH-1015 Lausanne, Switzerland.\\
{\tt jacques.thevenaz@epfl.ch}


\begin{thebibliography}{}

\bibitem[BT1]{BT1}
S.~Bouc, J.~Th\'evenaz.
\newblock The algebra of essential relations on a finite set, 
\newblock {\em J. reine angew. Math.} 712 (2016), 225--250.

\bibitem[BT2]{BT2}
S.~Bouc, J.~Th\'evenaz.
\newblock Correspondence functors and finiteness conditions,
\newblock {\em J. Algebra} 495 (2018), 150--198.

\bibitem[BT3]{BT3}
S.~Bouc, J.~Th\'evenaz.
\newblock Correspondence functors and lattices,
\newblock {\em J. Algebra} 518 (2019), 453--518.

\bibitem[BT4]{BT4}
S.~Bouc, J.~Th\'evenaz.
\newblock The algebra of Boolean matrices, correspondence functors, and simplicity,
\newblock preprint, 2018.

\bibitem[Ro]{Ro}
S.~Roman.
\newblock {\em Lattices and ordered sets},
\newblock Springer, New York, 2008.

\bibitem[St]{St}
R.~P.~Stanley.
\newblock {\em Enumerative Combinatorics,  Vol.~I}, Second edition,
\newblock Cambridge studies in advanced mathematics 49,
\newblock Cambridge University Press, 2012.


\end{thebibliography}
\end{document}